\documentclass[12pt]{article}
\usepackage{latexsym, amssymb, amsmath, amscd, amsfonts, epsfig, graphicx, colordvi,verbatim,ifpdf,extarrows}
\usepackage{amsfonts, amsmath, amssymb}
\usepackage{amssymb,amsfonts,amsmath,latexsym,epsfig,cite, psfrag,eepic,color}
\usepackage{amscd,graphics}
\usepackage{latexsym, amssymb,  amsmath,amscd, amsfonts, epsfig, graphicx, colordvi,amsthm}

\usepackage{graphicx}
\usepackage{epstopdf}
\usepackage{color}
\usepackage{ifpdf}
\usepackage{fancybox}
\usepackage[font=small,labelfont=bf,labelsep=none]{caption}
\usepackage{float}

\newtheorem{thm}{Theorem}[section]
\newtheorem{pro}[thm]{Proposition}

\newtheorem{defi}[thm]{Definition}
\newtheorem{lem}[thm]{Lemma}

\def\pf{\noindent{\it Proof.} }
\def\qed{\nopagebreak\hfill{\rule{4pt}{7pt}}
\medbreak}

\numberwithin{equation}{section}

\def\qed{\nopagebreak\hfill{\rule{4pt}{7pt}}
\medbreak}

\newlength{\boxedparwidth}
  {\begin{center} \begin{tabular}{|@{\hspace{.315in}}c@{\hspace{.15in}}|}
                  \hline \\ \begin{minipage}[t]{\boxedparwidth}
                  \setlength{\parindent}{.25in}}%
  {\end{minipage} \\ \\ \hline \end{tabular} \end{center}}

\begin{document}

\begin{center}

 {\Large \bf The Bressoud-G\"ollnitz-Gordon Theorem for \\[5pt] Overpartitions of even moduli}
\end{center}

\begin{center}
 {Thomas Y. He}$^{1}$, {Allison Y.F. Wang}$^{2}$ and
  {Alice X.H. Zhao}$^{3}$ \vskip 2mm

   $^{1,2}$Center for Applied Mathematics \\ Tianjin University, Tianjin 300072, P.R. China\\[6pt]

   $^{3}$Center for Combinatorics, LPMC-TJKLC\\
   Nankai University, Tianjin 300071, P.R. China\\[6pt]

   \vskip 2mm

    $^1$heyao@tju.edu.cn, $^2$wangyifang@tju.edu.cn, $^3$zhaoxiaohua@mail.nankai.edu.cn
\end{center}

\vskip 6mm \noindent {\bf Abstract.}  We give an overpartition analogue of Bressoud's combinatorial generalization of the G\"ollnitz-Gordon theorem for even moduli in general case. Let $\widetilde{O}_{k,i}(n)$ be the number of overpartitions of $n$ whose parts satisfy certain difference condition and $\widetilde{P}_{k,i}(n)$ be the number of overpartitions of $n$ whose non-overlined parts satisfy certain congruence condition. We show that $\widetilde{O}_{k,i}(n)=\widetilde{P}_{k,i}(n)$ for $1\leq i<k$.



\section{Introduction}
\ \ \ \ Throughout this article, we shall adopt the common notation as used in Andrews \cite{Andrews-1976}. Let
\[(a)_\infty=(a;q)_\infty=\prod_{i=0}^{\infty}(1-aq^i),\]
and
\[(a)_n=(a;q)_n=\frac{(a)_\infty}{(aq^n)_\infty}.\]

The Rogers-Ramanujan identities \cite{Rogers-1894} are the most fascinating identities in the theory of partitions. They can be stated either analytically or combinatorially. Each viewpoint has led to its own generalizations. On the combinatorial side, Gordon \cite{Gordon-1961} made the first break-through with an infinite family of identities in 1961 by proving a combinatorial generalization of the Rogers-Ramanujan identities, which is stated as the following theorem.

\begin{thm}[Rogers-Ramanujan-Gordon]\label{R-R-G}
For $k\geq i\geq1$. Let $B_{k,i}(n)$ denote the number of partitions $\lambda$ of $n$ of the form $\lambda=(1^{f_1}, 2^{f_{2}}, 3^{f_3},\ldots)$ where $f_t(\lambda)$ $($$f_t$ for short$)$ denotes the number of times the number $t$ appears as a part in $\lambda$ such that \begin{itemize}
\item[{\rm(1)}] $f_1(\lambda)\leq i-1$;
\item[{\rm (2)}] $f_t(\lambda)+f_{t+1}(\lambda)\leq k-1$.
\end{itemize}
Let $A_{k,i}(n)$ denote the number of partitions of $n$ into parts $\not\equiv0,\pm i\pmod{2k+1}$.

Then, for all $n\geq0$,
\[A_{k,i}(n)=B_{k,i}(n).\]
\end{thm}
Subsequently, Andrews \cite{Andrews-1974} discovered the generating function version of Theorem \ref{R-R-G}.

\begin{thm}[Andrews]\label{R-R-G-e}
For $k\geq i\geq 1$,
\begin{equation}\label{R-R-A1}
\sum_{N_1\geq N_2\geq \ldots \geq N_{k-1}\geq0}\frac{q^{N_1^2+N_2^2+\cdots+N^2_{k-1}+N_i+\cdots+N_{k-1}}}{(q)_{N_1-N_2}(q)_{N_2-N_3}\cdots(q)_{N_{k-1}}}=\frac{(q^i,q^{2k+1-i},q^{2k+1};q^{2k+1})_\infty}{(q)_\infty}.
\end{equation}
\end{thm}

By using the $q$-difference method, Andrews \cite{Andrews-1974} showed that the generating function of $B_{k,i}(n)$  defined in Theorem \ref{R-R-G} equals the left hand side of \eqref{R-R-A1}. More precisely,  he obtained the following formula for the generating function of   $B_{k,i}(m,n)$, where  $B_{k,i}(m,n)$ denotes the number of partitions enumerated by $B_{k,i}(n)$ that have $m$ parts. Recently, Kur\c{s}ung\"{o}z \cite{Kursungoz-2009,Kursungoz-2010} gave the generating function of $B_{k,i}(m,n)$ by introducing the notion of the Gordon marking of a partition.
\begin{thm}\label{R-R-G-left-e}
For $k\geq i\geq 1$,
\begin{equation}\label{R-R-A1-gf}
\sum_{n\geq 0}B_{k,i}(m,n)x^mq^n=\sum_{N_1\geq N_2\geq \ldots \geq N_{k-1}\geq0}\frac{q^{N_1^2+N_2^2+\cdots+N^2_{k-1}+N_i
+\cdots+N_{k-1}}x^{N_1+\cdots+N_{k-1}}}{(q)_{N_1-N_2}
(q)_{N_2-N_3}\cdots(q)_{N_{k-1}}}.
\end{equation}
\end{thm}

In 1979, Bressoud \cite{Bressoud-1979} extended Rogers-Ramanujan-Gordon theorem to even moduli.
\begin{thm}[Bressoud-Rogers-Ramanujan-Gordon] \label{R-R-B}
For  $k>i\geq 1$.  Let $\widetilde{A}_{k,i}(n)$ denote the number of partitions of $n$ of the form $\lambda=(1^{f_1}, 2^{f_{2}}, 3^{f_3},\ldots)$ such that
\begin{itemize}
\item[{\rm (1)}] $f_1(\lambda)\leq i-1$;

\item[{\rm (2)}] $f_t(\lambda)+f_{t+1}(\lambda)\leq k-1$;

\item[{\rm (3)}] If $f_t(\lambda)+f_{t+1}(\lambda)= k-1$, then $tf_t(\lambda)+(t+1)f_{t+1}(\lambda)\equiv i-1\pmod{2}$.
\end{itemize}
 Let $\widetilde{B}_{k,i}(n)$ denote the number of partitions of $n$ into parts $\not\equiv0,\pm i\pmod{2k}$.

 Then, for all $n\geq0$,
\[\widetilde{A}_{k,i}(n)=\widetilde{B}_{k,i}(n).\]
\end{thm}

In \cite{Bressoud-1980},  Bressoud found the following generating function version of Theorem \ref{R-R-B}.

\begin{thm}[Bressoud]\label{R-R-B-e}
 For $k>i\geq 1$.
\begin{equation}\label{R-R-BB1}
\sum_{N_1\geq N_2\geq \ldots \geq N_{k-1}\geq0}\frac{q^{N_1^2+N_2^2+\cdots+N^2_{k-1}+N_i+\cdots+N_{k-1}}}
{(q)_{N_1-N_2}(q)_{N_2-N_3}\cdots(q^{2};q^{2})_{N_{k-1}}}=\frac{(q^i,q^{2k-i},q^{2k};q^{2k})_\infty}{(q)_\infty}.
\end{equation}
\end{thm}

Let $\widetilde{B}_{k,i}(m,n)$ denote the number of partitions enumerated by $\widetilde{B}_{k,i}(n)$ that have $m$ parts. Using the notion of the Gordon marking of a partition, Kur\c{s}ung\"{o}z \cite{Kursungoz-2009,Kursungoz-2010} also obtained the generating function of $\widetilde{B}_{k,i}(m,n)$.
\begin{thm}\label{R-R-G-left-o}
For $k\geq i\geq 1$,
\begin{equation}\label{R-R-A1-gf-ee}
\sum_{n\geq 0}\widetilde{B}_{k,i}(m,n)x^mq^n=\sum_{N_1\geq N_2\geq \ldots \geq N_{k-1}\geq0}\frac{q^{N_1^2+N_2^2+\cdots+N^2_{k-1}+N_i
+\cdots+N_{k-1}}x^{N_1+\cdots+N_{k-1}}}{(q)_{N_1-N_2}
(q)_{N_2-N_3}\cdots(q^2;q^2)_{N_{k-1}}}.
\end{equation}
\end{thm}

G\"ollnitz-Gordon identities are independently introduced by Gordon\cite{Gordon-1962, Gordon-1965} and  G\"ollnitz\cite{Gollnitz-1960, Gollnitz-1967}. In 1967,   Andrews\cite{Andrews-1967}  found the combinatorial generalization of the G\"ollnitz-Gordon identities. In \cite{Bressoud-1980},  Bressoud extended Andrews-G\"ollnitz-Gordon theorem to even moduli.

\begin{thm}[Andrews-Bressoud-G\"ollnitz-Gordon]\label{Gordon-Gollnitz-even}
For $k>i\geq1$. Let $\widetilde{C}_{k,i}(n)$ denote the number of partitions of $n$ of the form $\lambda=(1^{f_1}, 2^{f_{2}}, 3^{f_3},\ldots)$ such that
\begin{itemize}
\item[{\rm (1)}] $f_1(\lambda)+f_2(\lambda)\leq i-1$;

\item[{\rm (2)}] $f_{2t+1}(\lambda)\leq1$;

\item[{\rm (3)}] $f_{2t}(\lambda)+f_{2t+1}(\lambda)+f_{2t+2}(\lambda)\leq k-1$;

\item[{\rm (4)}] if $f_{2t}(\lambda)+f_{2t+1}(\lambda)+f_{2t+2}(\lambda)=k-1$, then $tf_{2t}(\lambda)+tf_{2t+1}(\lambda)+(t+1)f_{2t+2}(\lambda)\equiv {O}_{2t+1}+i-1\pmod{2}$, where ${O}_{2t+1}(\lambda)$ denotes the number of odd parts in $\lambda$ those not exceed $2t+1$.
\end{itemize}

\noindent Let $\widetilde{D}_{k,i}(n)$ denote the number of partitions of $n$ into parts $\not\equiv2\pmod4$ and $\not\equiv0,\pm(2i-1)\pmod{4k-2}$.
Then, for all $n\geq0$,
\[\widetilde{C}_{k,i}(n)=\widetilde{D}_{k,i}(n).\]
\end{thm}

 Bressoud \cite{Bressoud-1980}
obtained the following generating function version of Theorem \ref{Gordon-Gollnitz-even}.

\begin{thm}[Bressoud]\label{Gollnitz-3-e}
For $k> i\geq1$.
\begin{equation}\label{Gollnitz-A2}
\begin{split}
&\sum_{N_1\geq N_2\geq \ldots \geq N_{k-1}\geq0}\frac{(-q^{1-2N_1};q^2)_{N_1}
q^{2(N^2_1+\cdots+N^2_{k-1}+N_i+\cdots+N_{k-1})}}
{(q^2;q^2)_{N_1-N_2}(q^2;q^2)_{N_2-N_3}\cdots(q^{4};q^{4})_{N_{k-1}}}\\
&=\frac{(q^2;q^4)_\infty(q^{4k-2},q^{2i-1},q^{4k-2i-1};q^{4k-2})_\infty}
{(q;q)_\infty}.
\end{split}
\end{equation}
\end{thm}

Recall that an overpartition $\lambda$ of  $n$ is a partition of $n$ in which the first occurrence of a number can be overlined.  In this paper, we write an overpartition $\lambda$ as the form $(a^{f_a},\bar{a}^{f_{\bar a}},\ldots, 1^{f_1},\bar{1}^{f_{\bar 1}})$  where $f_t$ (resp. $f_{\bar{t}}$) denotes the number of times the number $t$ (resp. $\bar{t}$)  appears as a part (resp. overlined part) in $\lambda$. In recent years, there are many overpartition analogues of classical partition theorems.  For example, Corteel and Mallet \cite{Corteel-2007}, Corteel, Lovejoy and Mallet \cite{Corteel-2008}, Lovejoy \cite{Lovejoy-2003, Lovejoy-2004, Lovejoy-2007,Lovejoy-2010} found many overpartition analogues of the Rogers-Ramanujan-Gordon theorem. Recently, Chen, Sang and Shi \cite{Chen-Sang-Shi-2013} obtained  overpartition analogue of  Rogers-Ramanujan-Gordon theorem in the general case.
 He, Ji, Wang and Zhao \cite{He-Ji-Wang-Zhao-2016} obtained the overpartition analogue of Andrews-G\"ollnitz-Gordon theorem.

Using the $q$-difference method, Chen, Sang and Shi \cite{Chen-Sang-Shi-2015} also obtained  overpartition analogue of Theorem \ref{R-R-B} in the general case.
\begin{thm}\label{BRESSOUD-even-overp}
For $k>i\geq1$. Let $\widetilde{E}_{k,i}(n)$ denote the number of overpartitions of $n$ of the form $(\bar{1}^{f_{\bar 1}},1^{f_1},\bar{2}^{f_{\bar 2}},2^{f_2},\ldots)$ such that
\begin{itemize}
\item[{\rm (1)}] $f_{{1}}(\lambda)\leq i-1$;

\item[{\rm (2)}] $f_{t}(\lambda)+f_{\overline{t}}(\lambda)+f_{{t+1}}(\lambda)\leq k-1$;

\item[{\rm (3)}] If $f_{t}(\lambda)+f_{\overline{t}}(\lambda)+f_{{t+1}}(\lambda)= k-1$ then $tf_{t}(\lambda)+tf_{\overline{t}}(\lambda)+(t+1)f_{t+1}(\lambda)\equiv V_{t}(\lambda)+i-1\pmod{2}$,
\end{itemize}
where $V_{s}(\lambda)$ denotes the number of overlined parts those not exceed $s$.

Let $\widetilde{F}_{k,i}(n)$ denote the number of overpartitions of $n$ whose non-overlined parts are not congruent to $0,\pm i$ modulo $2k-1$.

Then, for all $n\geq0$,
\[\widetilde{E}_{k,i}(n)=\widetilde{F}_{k,i}(n).\]
 \end{thm}

Chen, Sang and Shi \cite{Chen-Sang-Shi-2015} also gave generating function version of Theorem \ref{BRESSOUD-even-overp}.
 \begin{thm}\label{R-R-B-e-e-1}
 For $k>i\geq 1$.
\begin{equation}\label{R-R-BB1-e-1}
\begin{split}
&\sum_{N_1\geq N_2\geq \ldots \geq N_{k-1}\geq0}\frac{(-q^{1-N_1})_{N_1-1}q^{N_1^2+N_2^2+\cdots+N^2_{k-1}+N_{i+1}+\cdots+N_{k-1}}(1+q^{N_i})}
{(q)_{N_1-N_2}\cdots(q)_{N_{k-2}-N_{k-1}}(q^{2};q^{2})_{N_{k-1}}}\\
&=\frac{(-q)_\infty(q^i,q^{2k-i-1},q^{2k-1};q^{2k-1})_\infty}{(q)_\infty}.
\end{split}
\end{equation}
\end{thm}

Let $\widetilde{E}_{k,i}(m,n)$ denote the number of overpartitions counted by $\widetilde{E}_{k,i}(n)$ with exactly $m$ parts, Sang and Shi \cite{Sang-Shi-2015} gave a combinatorial proof of the following  generating function  of $\widetilde{E}_{k,i}(m,n)$.
 \begin{thm}\label{GENERATING-nwe-2}
For $k> i\geq1$, we have
\begin{equation}\label{GENERATING-nwe2}
\begin{split}
&\displaystyle\sum_{m,n\geq0}\widetilde{E}_{k,i}(m,n)x^mq^n\\
&=\sum_{N_1\geq N_2\geq \ldots \geq N_{k-1}\geq0}\frac{(-q^{1-N_1})_{N_1-1}
q^{N^2_1+\cdots+N^2_{k-1}+N_{i+1}+\cdots+N_{k-1}}(1+q^{N_i})x^{N_1+\cdots+N_{k-1}}}
{(q)_{N_1-N_2}\cdots(q)_{N_{k-2}-N_{k-1}}(q^{2};q^{2})_{N_{k-1}}}.
\end{split}
\end{equation}
\end{thm}

In this article, we will give an another proof of Theorem \ref{R-R-B-e-e-1} in the next section by using Bailey's lemma and the change of base formula and a direct combinatorial proof of Theorem \ref{GENERATING-nwe-2} by introduce the notation of Gordon marking for an overpartition introduced by Chen, Sang, and Shi \cite{Chen-Sang-Shi-2013}.

The main result of this article is the following overpartition analogue of Theorem \ref{Gordon-Gollnitz-even}.
\begin{thm}\label{Gordon-Gollnitz-even-overp}
For $k>i\geq1$. Let $\widetilde{O}_{k,i}(n)$ denote the number of overpartitions of $n$ of the form $(a^{f_a},\bar{a}^{f_{\bar a}},\ldots, 1^{f_1},\bar{1}^{f_{\bar 1}})$ such that
\begin{itemize}
\item[{\rm (1)}] $f_{\overline{1}}(\lambda)+f_2(\lambda)\leq i-1$;

\item[{\rm (2)}] $f_{2t}(\lambda)+f_{\overline{2t}}(\lambda)+f_{\overline{2t+1}}(\lambda)+f_{2t+2}(\lambda)\leq k-1$;

\item[{\rm (3)}]If $f_{2t+1}(\lambda)\geq1$, then $f_{2t+2}(\lambda)\leq k-2$;

\item[{\rm (4)}] If $f_{2t}(\lambda)+f_{\overline{2t}}(\lambda)+f_{\overline{2t+1}}(\lambda)+f_{2t+2}(\lambda)= k-1$ then $tf_{2t}(\lambda)+tf_{\overline{2t}}(\lambda)+tf_{\overline{2t+1}}(\lambda)+(t+1)f_{2t+2}(\lambda)\equiv V_{2t+1}(\lambda)+i-1\pmod{2}$;

\item[{\rm (5)}]If $f_{2t+1}(\lambda)\geq1$ and  $f_{2t+2}(\lambda)= k-2$, then $t+(t+1)f_{2t+2}(\lambda)\equiv V_{2t+1}(\lambda)+i-1\pmod{2}$,
\end{itemize}
where $V_{s}(\lambda)$ denotes the number of overlined parts those not exceed $s$.

Let $\widetilde{P}_{k,i}(n)$ denote the number of overpartitions of $n$ whose non-overlined parts are not congruent to $0,\pm (2i-1)$ modulo $4k-4$.

Then, for all $n\geq0$,
\[\widetilde{O}_{k,i}(n)=\widetilde{P}_{k,i}(n).\]
 \end{thm}

It should be noted that if an overpartition $\lambda$ counted by $\widetilde{O}_{k,i}(n)$ does not contain overlined even parts and non-overlined odd part, and we change the overlined odd parts in $\lambda$ to non-overlined parts, then we get a partition enumerated by $\widetilde{C}_{k,i}(n)$. Hence we say that Theorem \ref{Gordon-Gollnitz-even-overp} is the overpartition analogue of Theorem \ref{Gordon-Gollnitz-even}.

The corresponding generating function version of Theorem \ref{Gordon-Gollnitz-even-overp} is given as the following theorem.

\begin{thm}\label{Gordon-Gollnitz-even-overp-eqn}
For $k> i\geq1$.
\begin{equation}\label{over4}
\begin{split}
&\sum_{N_{1}\geq \cdots\geq N_{k-1}\geq 0}\frac{(-q^{2-2N_1};q^2)_{N_1-1}(-q^{1-2N_1};q^2)_{N_1}q^{2(N^2_{1}+\cdots
+N^2_{k-1}+N_{i+1}+\cdots+N_{k-1})}(1+q^{2N_i})}
{(q^2;q^2)_{N_1-N_2}\cdots(q^2;q^2)_{N_{k-2}-N_{k-1}}
(q^{4};q^{4})_{N_{k-1}}}\\
&=\frac{(-q;q)_\infty(q^{2i-1},q^{4k-3-2i},q^{4k-4};q^{4k-4})_\infty}{(q;q)_\infty}.
\end{split}
\end{equation}

\end{thm}

We will first give an analytic proof of Theorem \ref{Gordon-Gollnitz-even-overp-eqn} in the next section by using Bailey's lemma and the change of base formula. We then use Theorem \ref{Gordon-Gollnitz-even-overp-eqn} to derive  Theorem \ref{Gordon-Gollnitz-even-overp}. To be more precisely,  let $\widetilde{O}_{k,i}(m,n)$ denote the number of overpartitions counted by $\widetilde{O}_{k,i}(n)$ with exactly $m$ parts, we shall give a combinatorial proof of the following   generating function  of $\widetilde{O}_{k,i}(m,n)$ by introducing the G\"ollnitz-Gordon marking \cite{He-Ji-Wang-Zhao-2016} of an overpartition.

\begin{thm}\label{Gollnitz-even-e1}
For $k>i\geq1$, we have
\begin{equation}\label{Gollnitz-evene1}
\begin{split}
& \displaystyle\sum_{m,n\geq0}\widetilde{O}_{k,i}(m,n)x^mq^n\\
&=\sum_{N_{1}\geq \cdots\geq N_{k-1}\geq 0}\frac{(-q^{2-2N_1};q^2)_{N_1-1}(-q^{1-2N_1};q^2)_{N_1}
q^{2(N^2_1+\cdots+N^2_{k-1}+N_{i+1}+\cdots+N_{k-1})}
(1+q^{2N_i})x^{N_1+\cdots+N_{k-1}}}
{(q^2;q^2)_{N_1-N_2}\cdots(q^2;q^2)_{N_{k-2}-N_{k-1}}(q^{4};q^{4})_{N_{k-1}}}.\\
\end{split}
\end{equation}
\end{thm}

By setting $x=1$ in \eqref{Gollnitz-evene1}, we obtain the generating function for $\widetilde{O}_{k,i}(n)$ which is the left hand side of  \eqref{over4}. On the other hand, it is evident that the generating function of $\widetilde{P}_{k,i}(n)$ equals
\begin{equation}\label{Gollnitz-even-left}
\displaystyle\sum_{n\geq0}\widetilde{P}_{k,i}(n)q^n=\frac{(-q;q)_\infty(q^{2i-1},q^{4k-3-2i},q^{4k-4};q^{4k-4})_\infty}
{(q;q)_\infty},
\end{equation}
which is the right hand side of  \eqref{over4}. Hence we are led to  Theorem \ref{Gordon-Gollnitz-even-overp} by Theorem  \ref{Gordon-Gollnitz-even-overp-eqn}.

 This paper is organized as follows. In Section 2,  we show Theorem \ref{R-R-B-e-e-1} and Theorem \ref{Gordon-Gollnitz-even-overp-eqn} by using Bailey's lemma and the change of base formula due to Bressoud, Ismail and Stanton. In Section 3, we recall the definition and related notations of Gordon marking of an overpartition and give the outline of the proof of Theorem \ref{GENERATING-nwe-2}. In Section 4, we show the detail of the proof of Theorem \ref{GENERATING-nwe-2}. In Section 5, we recall the definition and related notations of G\"ollnitz-Gordon marking of an overpartition and give an outline of the proof of Theorem \ref{Gollnitz-even-e1}. In Section 6, we show the proof of Theorem \ref{Gollnitz-even-e1}, thus we complete the proof of Theorem \ref{Gordon-Gollnitz-even-overp}.

\section{Proof of Theorem \ref{R-R-B-e-e-1} and Theorem \ref{Gordon-Gollnitz-even-overp-eqn}}

 We first briefly review Bailey pairs and   Bailey's lemma. Recall that a pair of sequences $(\alpha_n(a,q),\beta_n(a,q))$ is called a Bailey pair with parameters $(a,q)$ if they have the following relation for all $n\geq 0,$
\begin{equation}\label{bailey pair}
\beta_n(a,q)=\sum_{r=0}^n\frac{\alpha_r(a,q)}{(q;q)_{n-r}(aq;q)_{n+r}}.
\end{equation}

 Bailey's lemma was first given by Bailey \cite{Bailey-1949}  and was formulated by Andrews  \cite{Andrews-1984, Andrews-1986} in the following form.

\begin{thm}[Bailey's lemma]\label{bailey lemma}
If $(\alpha_n(a,q),\beta_n(a,q))$ is a Bailey pair with parameters $(a,q)$,
then $(\alpha'_n(a,q),\beta'_n(a,q))$ is another Bailey pair with parameters $(a,q)$, where
\begin{equation}\label{bl}
\begin{split}
\alpha_n'(a,q)&=\frac{(\rho_1;q)_n(\rho_2;q)_n}{(aq/\rho_1;q)_n(aq/\rho_2;q)_n}\left(\frac{aq}
{\rho_1\rho_2}\right)^n\alpha_n(a,q),\\[5pt]
\beta_n'(a,q)&=\sum_{j=0}^{n}\frac{(\rho_1;q)_j(\rho_2;q)_j(aq/\rho_1\rho_2;q)_{n-j}}
{(aq/\rho_1;q)_n(aq/\rho_2;q)_n(q;q)_{n-j}}\left(\frac{aq}{\rho_1\rho_2}\right)^j\beta_j(a,q).
\end{split}
\end{equation}
\end{thm}
Andrews first noticed that Bailey's lemma can create a new Bailey pair from a given one. Hence the iteration of the lemma leads to a sequence of Bailey pairs called a Bailey chain. Based on this observation,  Andrews \cite{Andrews-1984} used Bailey's lemma to
show the Andrews-Gordon identity \eqref{R-R-A1} in Theorem \ref{R-R-G-e} and the Bressoud-Gordon identity \eqref{R-R-BB1} in Theorem \ref{R-R-B-e} with $i=1$ and $i=k$. Subsequently, Agarwal, Andrews and Bressoud \cite{Agarwal- Andrews-Bressoud} gave an extension of the Bailey chain known as the Bailey lattice, by means of which enable us to prove  the Andrews-Gordon identity and the Bressoud-Gordon identity \eqref{R-R-BB1} with $1\leq i\leq k$. In \cite{Bressoud-Ismail-Stanton-2000},  Bressoud, Ismail and Stanton  established versions of Bailey's lemma, known as the change of base formulas,  which can be used  to  prove the Bressound-G\"ollnitz-Gordon identity \eqref{Gollnitz-A2}.

In this section,  we will show Theorem \ref{R-R-B-e-e-1} and Theorem \ref{Gordon-Gollnitz-even-overp-eqn} by combining Bailey's lemma and the change of base formula. First, we need to review two limiting cases of Bailey's lemma which are both appeared in \cite{Paule-1985,Paule-1987}. The first limiting case is obtained by letting $\rho_1,\rho_2\rightarrow \infty$ in Theorem \ref{bailey lemma}.

\begin{lem}\label{bl1}
If $(\alpha_n(a,q),\beta_n(a,q))$ is a Bailey pair with parameters $(a,q)$,
then $(\alpha_n', \beta_n')$ is also a Bailey pair with parameters $(a,q)$, where
\begin{equation}
\begin{split}
&\alpha_n'(a,q)=a^nq^{n^2}\alpha_n(a,q),\\[3pt]
&\beta_n'(a,q)=\sum_{j=0}^n\frac{a^jq^{j^2}}{(q)_{n-j}}\beta_j(a,q).
\end{split}
\end{equation}
\end{lem}

 The second limiting case is obtained by letting $\rho_1\rightarrow \infty,\rho_2=-1$ in Theorem \ref{bailey lemma}.

\begin{lem}\label{bl11}
If $(\alpha_n(a,q),\beta_n(a,q))$ is a Bailey pair with parameters $(a,q)$,
then $(\alpha_n', \beta_n')$ is also a Bailey pair with parameters $(a,q)$, where
\begin{equation}
\begin{split}
&\alpha_n'(a,q)=\frac{(-1)_na^nq^{(n^2+n)/2}}{(-aq)_n}\alpha_n(a,q),\\[3pt]
&\beta_n'(a,q)=\frac{1}{(-aq)_n}\sum_{j=0}^n\frac{(-1)_ja^jq^{(j^2+j)/2}}{(q)_{n-j}}\beta_j(a,q).
\end{split}
\end{equation}
\end{lem}

The proof of  Theorem \ref{Gordon-Gollnitz-even-overp-eqn} require  the following  special case of change of base formula.

\begin{lem}\label{bailey lemma3}{\rm \cite[Theorem 2.5, $B\rightarrow \infty$]{Bressoud-Ismail-Stanton-2000}}
If $(\alpha_n(a,q),\beta_n(a,q))$ is a Bailey pair with parameters $(a,q)$,
then $(\alpha_n', \beta_n')$ is also a Bailey pair with parameters $(a,q)$, where
\begin{equation*}
\begin{split}
  \alpha_n'(a,q)&=\frac{1+a}{1+aq^{2n}}q^n\alpha_n(a^2,q^2),\\[3pt]
  \beta_n'(a,q)&=\sum_{k=0}^n\frac{(-a;q)_{2k}q^k}{(q^2;q^2)_{n-k}}\beta_{k}(a^2,q^2).
\end{split}
\end{equation*}
\end{lem}

The following proposition  is  useful in our proofs as well.

\begin{pro}\label{bl4}{\rm \cite[Proposition 4.1]{Bressoud-Ismail-Stanton-2000}}
If $(\alpha_n,\beta_n)$ is a Bailey pair with parameters $(1,q)$, and
\[\alpha_n=\left\{
             \begin{array}{ll}
               1, & \hbox{for $n=0$}, \\[5pt]
               (-1)^nq^{An^2}(q^{(A-1)n}+q^{-(A-1)n}), & \hbox{for $n>0$,}
             \end{array}
           \right.
\]
then $(\alpha_n'(q),\beta_n'(q))$ is a Bailey pair with parameters $(1,q)$, where
\begin{align*}
 \alpha_n'(q)&=\left\{
             \begin{array}{ll}
               1, & \hbox{for $n=0$}, \\[3pt]
               (-1)^nq^{An^2}(q^{An}+q^{-An}), & \hbox{for $n>0$,}
             \end{array}
           \right. \\[10pt]
 \beta_n'(q)&=q^n\beta_n(q).
\end{align*}
\end{pro}

\noindent{\bf Proof of Theorem \ref{R-R-B-e-e-1}: } We begin with the following  Bailey pair \cite[E(4)]{Slater-1952}  with parameters $(1,q)$, where
\begin{equation*}
\begin{split}
&\alpha^{(1)}_n(1,q)=\left\{
             \begin{array}{ll}
               1, & \hbox{if $n=0$}, \\[5pt]
               (-1)^nq^{n^2}(q^{-n}+q^{n}), & \hbox{if $n\geq1$,}
             \end{array}
               \right.\\[5pt]
&\beta^{(1)}_n(1,q)=\frac{q^n}{(q^2;q^2)_n}.
\end{split}
\end{equation*}

Invoking Lemma \ref{bl1} gives the following Bailey pair
\begin{equation*}
\begin{split}
&\alpha^{(2)}_n(1,q)=\left\{
             \begin{array}{ll}
               1, & \hbox{if $n=0$}, \\[5pt]
               (-1)^nq^{2n^2}(q^{-n}+q^{n}), & \hbox{if $n\geq1$,}
             \end{array}
               \right.\\[10pt]
&\beta^{(2)}_n(1,q)=\sum_{N_{k-1}=0}^n\frac{q^{N_{k-1}^2+N_{k-1}}}
{(q)_{n-N_{k-1}}(q^2;q^2)_{N_{k-1}}}.
\end{split}
\end{equation*}

Alternately apply Proposition \ref{bl4} and Lemma \ref{bl1} for $k-i-2$ times. The result is
\begin{equation}\label{BP7}
\begin{split}
&\alpha^{(2k-2i-2)}_n(1,q)=\left\{
             \begin{array}{ll}
               1, & \hbox{if $n=0$}, \\[5pt]
               (-1)^nq^{(k-i)n^2}(q^{-(k-i-1)n}+q^{(k-i-1)n}), & \hbox{if $n\geq1$,}
             \end{array}
               \right.\\[10pt]
&\beta^{(2k-2i-2)}_n(1,q)=\sum_{n\geq N_{i+1}\geq \cdots \geq N_{k-1}\geq0}\frac{q^{N_{i+1}^2+\cdots+N_{k-1}^2+N_{i+1}+\cdots+N_{k-1}}}
{(q)_{n-N_{i+1}}(q)_{N_{i+1}-N_{i+2}}\cdots(q)_{N_{k-2}-N_{k-1}}(q^2;q^2)_{N_{k-1}}}.
\end{split}
\end{equation}

Applying Proposition \ref{bl4} to the Bailey pair \eqref{BP7} yields
\begin{equation}\label{BP8}
\begin{split}
&\alpha^{(2k-2i-1)}_n(1,q)=\left\{
             \begin{array}{ll}
               1, & \hbox{if $n=0$}, \\[5pt]
               (-1)^nq^{(k-i)n^2}(q^{-(k-i)n}+q^{(k-i)n}), & \hbox{if $n\geq1$,}
             \end{array}
               \right.\\[10pt]
&\beta^{(2k-2i-1)}_n(1,q)=q^n\sum_{n\geq N_{i+1}\geq \cdots \geq N_{k-1}\geq0}\frac{q^{N_{i+1}^2+\cdots+N_{k-1}^2+N_{i+1}+\cdots+N_{k-1}}}
{(q)_{n-N_{i+1}}(q)_{N_{i+1}-N_{i+2}}\cdots(q)_{N_{k-2}-N_{k-1}}(q^2;q^2)_{N_{k-1}}}.
\end{split}
\end{equation}

Adding these two Bailey pairs \eqref{BP7} and \eqref{BP8} together according to the definition of Bailey pair, we get a new Bailey pair relative to $(1,q)$
\begin{equation}\label{BP9}
\begin{split}
&\alpha^{(2k-2i)}_n(1,q)=\left\{
             \begin{array}{ll}
               1, & \hbox{if $n=0$}, \\[5pt]
               (-1)^nq^{(k-i)n^2}(q^{-(k-i)n}+q^{(k-i-1)n})(1+q^n)/2, & \hbox{if $n\geq1$,}
             \end{array}
               \right.\\[10pt]
&\beta^{(2k-2i)}_n(1,q)=\sum_{n\geq N_{i+1}\geq \cdots \geq N_{k-1}\geq0}\frac{(1+q^n)q^{N_{i+1}^2+\cdots+N_{k-1}^2+N_{i+1}+\cdots+N_{k-1}}}
{2(q)_{n-N_{i+1}}(q)_{N_{i+1}-N_{i+2}}\cdots(q)_{N_{k-2}-N_{k-1}}(q^2;q^2)_{N_{k-1}}}.
\end{split}
\end{equation}

Then apply Lemma \ref{bl1} to \eqref{BP9} $i-1$ times to get a Bailey pair $(\alpha^{(2k-i-1)}_n(1,q),\beta^{(2k-i-1)}_n(1,q))$ with parameters $(1,q)$, where
\begin{equation}\label{BP10}
\begin{split}
&\alpha^{(2k-i-1)}_n(1,q)=\left\{
             \begin{array}{ll}
               1, & \hbox{if $n=0$}, \\[5pt]
               (-1)^nq^{(k-1)n^2}(q^{-(k-i)n}+q^{(k-i-1)n})(1+q^n)/2 & \hbox{if $n\geq1$,}
             \end{array}
               \right.\\[10pt]
&\beta^{(2k-i-1)}_n(1,q)=\sum_{n\geq N_{2}\geq \cdots\geq N_{k-1}\geq 0}\frac{q^{N^2_{2}+N^2_{3}+\cdots+N^2_{k-1}+N_{i+1}+\cdots+N_{k-1}}(1+q^{N_i})}
{2(q)_{n-N_{2}}(q)_{N_{2}-N_{3}}\cdots(q)_{N_{k-2}-N_{k-1}}(q^2;q^2)_{N_{k-1}}}.
\end{split}
\end{equation}

Invoking Lemma \ref{bl11} to \eqref{BP10} with parameters $(1,q)$, gives
\begin{equation*}
\begin{split}
&\alpha^{(2k-i)}_n(1,q)=\left\{
             \begin{array}{ll}
               1, & \hbox{if $n=0$}, \\[5pt]
               (-1)^nq^{(2k-1)n^2/2}(q^{-(2k-2i-1)n/2}+q^{(2k-2i-1)n/2}), & \hbox{if $n\geq1$,}
             \end{array}
               \right.\\[10pt]
&\beta^{(2k-i)}_n(1,q)=\frac{1}{(-q)_n}\sum_{n\geq N_{1}\geq \cdots\geq N_{k-1}\geq 0}\frac{(-1)_{N_1}q^{(N_1+1)N_1/2+N^2_{2}+N^2_{3}+\cdots
+N^2_{k-1}+N_{i+1}+\cdots+N_{k-1}}(1+q^{N_i})}
{2(q)_{n-N_1}(q)_{N_1-N_{2}}\cdots
(q)_{N_{k-2}-N_{k-1}}(q^2;q^2)_{N_{k-1}}}.
\end{split}
\end{equation*}

By the definition of Bailey pairs, we have
\begin{equation}\label{biaohaohao}
\begin{split}
&\frac{1}{(-q)_n}\sum_{n\geq N_{1}\geq \cdots\geq N_{k-1}\geq 0}\frac{(-1)_{N_1}q^{(N_1+1)N_1/2+N^2_{2}+N^2_{3}+\cdots
+N^2_{k-1}+N_{i+1}+\cdots+N_{k-1}}(1+q^{N_i})}
{2(q)_{n-N_1}(q)_{N_1-N_{2}}\cdots
(q)_{N_{k-2}-N_{k-1}}(q^2;q^2)_{N_{k-1}}}\\
&=\frac{1}{(q)_n(q)_n}+\sum_{r=1}^{n}\frac{(-1)^rq^{(2k-1)r^2/2}(q^{-(2k-2i-1)r/2}+q^{(2k-2i-1)r/2})}{(q)_{n-r}(q)_{n+r}}
\end{split}
\end{equation}

Letting $n\rightarrow \infty$ and multiplying both sides by $(q^2;q^2)_\infty$ in \eqref{biaohaohao}, we obtain
\begin{equation}\label{zz3}
\begin{split}
&\sum_{N_{1}\geq \cdots\geq N_{k-1}\geq 0}\frac{(-q)_{N_1-1}q^{(N_1+1)N_1/2+N^2_{2}+N^2_{3}+\cdots
+N^2_{k-1}+N_{i+1}+\cdots+N_{k-1}}(1+q^{N_i})}
{(q)_{N_1-N_{2}}\cdots
(q)_{N_{k-2}-N_{k-1}}(q^2;q^2)_{N_{k-1}}}\\
&=\frac{(-q)_\infty}{(q)_\infty}\left(1+\sum_{n=1}^\infty
(-1)^nq^{(2k-1)n^2/2}(q^{-(2k-2i-1)n/2}+q^{(2k-2i-1)n/2})\right)
\end{split}
\end{equation}

\noindent Letting $q\rightarrow q^{(2k-1)/2},\ z\rightarrow-q^{(2k-2i-1)/2}$ in the following Jacobi's triple product identity.
\begin{equation}\label{Jacobi triple identity-1}
1+\sum_{n=1}^{\infty}q^{n^2}(z^{-n}+z^n)=(-zq;q^2)_\infty(-q/z;q^2)_\infty(q^2;q^2)_\infty.
\end{equation}
We derive that
\begin{equation}\label{zz1}
1+\sum_{n=1}^\infty
(-1)^nq^{(2k-1)n^2/2}(q^{-(2k-2i-1)n/2}+q^{(2k-2i-1)n/2})=(q^{i},q^{2k-i-1},q^{2k-1};q^{2k-1})_\infty.
\end{equation}

Submitting \eqref{zz1} into \eqref{zz3}, and noting that
\[(-q)_{N_1-1}q^{(N_1+1)N_1/2}=(-q^{1-N_1})_{N_1-1}q^{N_1^2},\]
we obtain \eqref{R-R-BB1-e-1}. Thus we complete the proof of Theorem \ref{R-R-B-e-e-1}. \qed

\noindent{\bf Proof of Theorem \ref{Gordon-Gollnitz-even-overp-eqn}: } Invoking Lemma \ref{bailey lemma3} to \eqref{BP10} with parameters $(1,q)$, we get another Bailey pair $(\alpha^{(2k-i)}_n(1,q),\beta^{(2k-i)}_n(1,q))$ with parameters $(1,q)$, where
\begin{equation*}
\begin{split}
&\alpha^{(2k-i)}_n(1,q)=\left\{
             \begin{array}{ll}
               1, & \hbox{if $n=0$}, \\[5pt]
               (-1)^nq^{2(k-1)n^2}(q^{-(2k-2i-1)n}+q^{(2k-2i-1)n}), & \hbox{if $n\geq1$,}
             \end{array}
               \right.\\[10pt]
&\beta^{(2k-i)}_n(1,q)=\sum_{n\geq N_{1}\geq \cdots\geq N_{k-1}\geq 0}\frac{(-q;q)_{2N_1-1}q^{N_1+2(N^2_{2}+N^2_{3}+\cdots
+N^2_{k-1}+N_{i+1}+\cdots+N_{k-1})}(1+q^{2N_i})}
{(q^2;q^2)_{n-N_1}(q^2;q^2)_{N_1-N_{2}}\cdots
(q^2;q^2)_{N_{k-2}-N_{k-1}}(q^4;q^4)_{N_{k-1}}}.
\end{split}
\end{equation*}

By the definition of Bailey pairs, letting $n\rightarrow \infty$ and multiplying both sides by $(q^2;q^2)_\infty$, we obtain
\begin{equation}\label{zz4}
\begin{split}
&\sum_{N_{1}\geq \cdots\geq N_{k-1}\geq 0}\frac{(-q)_{2N_1-1}q^{N_1+2(N^2_{2}+N^2_{3}+\cdots+N^2_{k-1}+N_{i+1}+\cdots+N_{k-1})}(1+q^{2N_i})}
{(q^2;q^2)_{N_1-N_{2}}\cdots(q^2;q^2)_{N_{k-2}-N_{k-1}}(q^4;q^4)_{N_{k-1}}}\\
&=\frac{(-q)_\infty}{(q)_\infty}\left(1+\sum_{n=1}^\infty
(-1)^nq^{2(k-1)n^2}(q^{-(2k-2i-1)n}+q^{(2k-2i-1)n})\right)
\end{split}
\end{equation}

Letting $q\rightarrow q^{2(k-1)},\ z\rightarrow-q^{(2k-2i-1)}$ in \eqref{Jacobi triple identity-1}, we derive that
\begin{equation}\label{zz2}
1+\sum_{n=1}^\infty
(-1)^nq^{2(k-1)n^2}(q^{-(2k-2i-1)n}+q^{(2k-2i-1)n})=(q^{2i-1},q^{4k-2i-3},q^{4k-4};q^{4k-4})_\infty.
\end{equation}
Submitting \eqref{zz2} into \eqref{zz4}, and noting that
\[(-q)_{2N_1-1}q^{N_1}=(-q^{2-2N_1};q^2)_{N_1-1}(-q^{1-2N_1;q^2})_{N_1}q^{2N^2_1},\]

We arrive at \eqref{over4}. Thus, we complete the proof of Theorem \ref{Gordon-Gollnitz-even-overp-eqn}. \qed

\section{Outline of proof of Theorem \ref{GENERATING-nwe-2}}

In this section, we give an outline of the proof of Theorem \ref{GENERATING-nwe-2}. Let $\mathbb{\widetilde{E}}_{k,i}(m,n)$ denote the set of overpartitions counted by $\widetilde{E}_{k,i}(m,n)$, we further classify $\mathbb{\widetilde{E}}_{k,i}(m,n)$ by considering whether the smallest part is overlined or not. Note that the parts of an overpartition are ordered in the following order.
 \begin{equation}\label{order}
 \overline{1}<1<\overline{2}<2<\cdots.
 \end{equation}
 Let $\mathbb{\widetilde{G}}_{k,i}(m,n)$ denote the set of overpartitions in $\mathbb{\widetilde{E}}_{k,i}(m,n)$ for which the smallest part is non-overlined, and let $\mathbb{\widetilde{H}}_{k,i}(m,n)$ denote the set of overpartitions in $\mathbb{\widetilde{E}}_{k,i}(m,n)$ with the smallest part is overlined.

Let $\widetilde{G}_{k,i}(m,n)=|\mathbb{\widetilde{G}}_{k,i}(m,n)|$ and $\widetilde{H}_{k,i}(m,n)=|\mathbb{\widetilde{H}}_{k,i}(m,n)|$. Sang and Shi obtained \cite{Sang-Shi-2015} a relation between ${\widetilde{G}}_{k,i}(m,n)$ and ${\widetilde{H}}_{k,i}(m,n)$.
\begin{thm}\label{2N1}
For $k> i\geq2$, we have
\begin{equation}\label{2N-3}
{\widetilde{G}}_{k,i}(m,n)={\widetilde{H}}_{k,i-1}(m,n).
\end{equation}
For $k>i=1$, we have
\begin{equation}\label{2N-4}
{\widetilde{G}}_{k,1}(m,n)={\widetilde{H}}_{k,k-1}(m,n-m).
\end{equation}
\end{thm}

{\pf}We sketch the proof of Theorem \ref{2N1} in \cite{Sang-Shi-2015} but omit the details here.

To prove \eqref{2N-4}, we give a bijection between $\mathbb{\widetilde{G}}_{k,1}(m,n)$ and $\mathbb{\widetilde{H}}_{k,k-1}(m,n-m)$. For an
overpartition $\lambda$ in $\mathbb{\widetilde{G}}_{k,1}(m,n)$, there are no parts equal to $1$ in $\lambda$, that is, each part is greater
than or equal to $2$, so we can substract $1$ from each part of ¦Ë and set one of the smallest
parts to an overlined part to obtain an overpartition $\lambda'$ in $\mathbb{\widetilde{H}}_{k,k-1}(m,n-m)$. Conversely, for an overpartition in $\mathbb{\widetilde{H}}_{k,k-1}(m,n-m)$, we can add $1$ to each part and change
the smallest overlined part to a non-overlined part to get an overpartition in $\mathbb{\widetilde{G}}_{k,1}(m,n)$. So
${\widetilde{G}}_{k,1}(m,n)={\widetilde{H}}_{k,k-1}(m,n-m)$ for all $k\geq2$.

For the case $k>i\geq2$, there is a simple bijection between $\mathbb{\widetilde{H}}_{k,i-1}(m,n)$ and $\mathbb{\widetilde{G}}_{k,i}(m,n)$.
Let $\lambda$ be an overpartition in $\mathbb{\widetilde{H}}_{k,i-1}(m,n)$. Switching the smallest overlined part of $\lambda$ to a
non-overlined part, we get an overpartition $\lambda'$  in $\mathbb{\widetilde{G}}_{k,i}(m,n)$. Conversely, for an overpartition in $\mathbb{\widetilde{G}}_{k,i}(m,n)$, we switch the smallest non-overlined part to an overlined part, we get an overpartition in $\mathbb{\widetilde{H}}_{k,i-1}(m,n)$. So
\eqref{2N-3} is established for $k>i\geq2$.    \qed

Recall that the Gordon marking of an overpartition $\lambda$ is an assignment of positive integers, called marks, to parts of $\lambda=(\lambda_1,\lambda_2,\ldots,\lambda_{\ell})$ where $1\leq\lambda_1\leq\lambda_2\leq\cdots\leq\lambda_{\ell}$, subject to certain conditions. More precisely, we assign the marks to parts also in the order stated in \eqref{order} such that the marks are as small as possible subject to the following conditions:

\begin{itemize}
     \item[{\rm (i).}] If $\overline{t+1}$ is not a part of $\lambda$, then all the parts $t,\overline{t}$, and $t+1$ are assigned different integers.
   \item[{\rm (ii).}] If $\lambda$ contains an overlined part $\overline{t+1}$, then the smallest mark assigned to a part $t$ or $\overline{t}$ can be used as the mark of $t+1$ or $\overline{t+1}$.
     \end{itemize}

For example, let
\[\lambda=(1,2,2,\bar{4},5,5,\bar{6},6,7,\bar{8},8,8,\overline{10},11,12,12,13,16).\]
Then the Gordon marking of $\lambda$ is
\[G(\lambda)=(1_1,2_2,2_3,\bar{4}_1,5_2,5_3,\bar{6}_1,6_2,7_3,\bar{8}_1,8_2,8_3,\overline{10}_1,{11}_2,{12}_1,{12}_3,{13}_2,{16}_1),\]
where the subscripts stand for marks. The Gordon marking of $\lambda$ can also be illustrated as follows.
\[G(\lambda)=\begin{array}{cc}\left[
\begin{array}{cccccccccccccccc}
&2&&&5&&7&8&&&&12&&&&\\
&2&&&5&6&&8&&&11&&13&&&\\
1&&&\bar{4}&&\bar{6}&&\bar{8}&&\overline{10}&&12&&&&16
\end{array}\right]&\begin{array}{c}3\\2\\1\end{array}
\end{array},\]
where column indicates the value of parts, and the row (counted from bottom to top) indicates the mark.

Let $N_r$ be the number of  parts in the $r$-th row of $G(\lambda)$. It is easy to find that $N_1\geq N_2\geq\cdots$. For the example above, we have $N_1=7,\ N_2=6,\ N_3=5$.

 Now, Let $\mathbb{\widetilde{E}}_{N_1,\ldots,N_{k-1};i}(n)$, $\mathbb{\widetilde{G}}_{N_1,\ldots,N_{k-1};i}(n)$ and $\mathbb{\widetilde{H}}_{N_1,\ldots,N_{k-1};i}(n)$  denote the sets of overpartitions $\lambda$ in $\mathbb{\widetilde{E}}_{k,i}(m,n)$, $\mathbb{\widetilde{G}}_{k,i}(m,n)$ and $\mathbb{\widetilde{H}}_{k,i}(m,n)$ respectively that have $N_r$ $r$-marked parts in their Gordon marking $G(\lambda)$ for $1\leq r\leq k-1$, where $\sum_{r=1}^{k-1}N_r=m$. Obviously, we have
 \begin{equation}\label{2N-1}
\mathbb{\widetilde{E}}_{N_1,\ldots,N_{k-1};i}(n)=\mathbb{\widetilde{G}}_{N_1,\ldots,N_{k-1};i}(n)\cup\mathbb{\widetilde{H}}_{N_1,\ldots,N_{k-1};i}(n).
\end{equation}

Let ${\widetilde{G}}_{N_1,\ldots,N_{k-1};i}(n)=|\mathbb{\widetilde{G}}_{N_1,\ldots,N_{k-1};i}(n)|$ and ${\widetilde{H}}_{N_1,\ldots,N_{k-1};i}(n)=|\mathbb{\widetilde{H}}_{N_1,\ldots,N_{k-1};i}(n)|$. Using the same method proving Theorem \ref{2N1}, we can easily get a relation between ${\widetilde{G}}_{N_1,\ldots,N_{k-1};i}(n)$ and ${\widetilde{H}}_{N_1,\ldots,N_{k-1};i}(n)$.
\begin{lem}\label{2Nnew1}
For $k> i\geq2$, we have
\begin{equation}\label{2Nnew-3}
{\widetilde{G}}_{N_1,\ldots,N_{k-1};i}(n)={\widetilde{H}}_{N_1,\ldots,N_{k-1};i-1}(n).
\end{equation}
For $k>i=1$, we have
\begin{equation}\label{2Nnew-4}
{\widetilde{G}}_{N_1,\ldots,N_{k-1};1}(n)={\widetilde{H}}_{N_1,\ldots,N_{k-1};k-1}(n-N_1-\cdots-N_{k-1}).
\end{equation}
\end{lem}

Let $\mathbb{\widetilde{B}}_{N_1,\ldots,N_{k-1};i}(n)$ denote the set of overpartitions in $\mathbb{\widetilde{G}}_{N_1,\ldots,N_{k-1};i}(n)$  for which there do not exist overlined part.

Set
\[\mathbb{\widetilde{G}}_{N_1,\ldots,N_{k-1};i}
=\bigcup_{n\geq0}\mathbb{\widetilde{G}}_{N_1,\ldots,N_{k-1};i}(n),\]
\[\mathbb{\widetilde{B}}_{N_1,\ldots,N_{k-1};i}
=\bigcup_{n\geq0}\mathbb{\widetilde{B}}_{N_1,\ldots,N_{k-1};i}(n).\]
We shall give a bijection for the following relation.

\begin{thm}\label{N-1}
For $N_1\geq N_2\geq \cdots \geq N_{k-1}\geq 0$, we have
\begin{equation}\label{N1}
\sum_{\lambda\in\mathbb{\widetilde{G}}_{N_1,\ldots,N_{k-1};i}}x^{\ell(\lambda)}
q^{|\lambda|}=(-q^{1-N_1})_{N_1-1}\sum_{\mu\in
\mathbb{\widetilde{B}}_{N_1,\ldots,N_{k-1};i}}
x^{\ell(\mu)}
q^{|\mu|},
\end{equation}
where $\ell(\lambda)$ denotes the number of parts of $\lambda$.
\end{thm}

\section{Proof of Theorem \ref{GENERATING-nwe-2}}

We modify the definitions of the first reduction operation and its reverse the first dilation operation which are introduced by Chen, Sang and Shi \cite{Chen-Sang-Shi-2013}. Here we just change the conditions that the first reduction operation and the first dilation operation must satisfy.

Let $\lambda=(\lambda_1,\ldots,\lambda_m)\in{\widetilde{G}}_{N_1,\ldots,N_{k-1};i}$ be an overpartition of $n$. we consider the Gordon marking of $\lambda$, then there are $N_1$ $1$-marked parts in $G(\lambda)$, set $1\leq\lambda_{N_1}^{(1)}\leq\cdots\leq\lambda_1^{(1)}$. If $\lambda_p^{(1)}$ is a non-overlined part and $\lambda_{p-1}^{(1)}$ is an overlined part for $1<p<N_1$ or $\lambda_p^{(1)}$ is a non-overlined part for $p=1$, we can define the first reduction operation of $p$-th. We assume $\lambda_p^{(1)}=t$, we consider the following two cases.

\noindent{\bf The first reduction operation of $p$-th kind:}

{\bf Case 1:} There is a non-overlined part $t+1$ of $\lambda$, but there is no overlined $1$-marked part $\overline{t+1}$. First, we change the part $\lambda_p^{(1)}$ to a $1$-marked overlined part $\bar{t}$.  Then, we
choose the part $t+1$ with the smallest mark, say $r$, and replace this $r$-marked part
$t+1$ with an $r$-marked part $t$. Moreover, if there is a $1$-marked overlined part to
the right of $\bar{t}$, then we switch it to a non-overlined part.

{\bf Case 2:} Either there is a $1$-marked overlined part $\overline{t+1}$ or there are no parts with
underlying part $t+1$. In this case, we may change the part $\lambda_p^{(1)}$ to a $1$-marked overlined part $\overline{t-1}$. Moreover, if there are $1$-marked parts larger than $t$, then we switch the overlined $1$-marked part next to $\lambda_p^{(1)}$ to a non-overlined part.

For example, let $\lambda$ be the following overpartition in ${\widetilde{G}}_{7,6,5;2}(135)$:
\[\lambda=\begin{array}{cc}\left[
\begin{array}{cccccccccccccccc}
&2&&&5&&7&8&&&&12&&&\\
&2&&&5&6&&8&&&11&&\mathbf{13}&&\\
1&&&\bar{4}&&\bar{6}&&\bar{8}&&\overline{10}&&\mathbf{12}&&&\mathbf{\overline{15}}
\end{array}\right]&\begin{array}{c}3\\2\\1\end{array}
\end{array},\]
It is easy to check that $\lambda^{(1)}_2=12$ is a non-overlined part and $\lambda^{(1)}_1=\overline{15}$ is an overlined part, so we can do The first reduction operation of $2$-th kind. Notice that $\overline{13}$ is not a $1$-marked part of $\lambda$, but $13$ is a $2$-marked part. By the operation in Case 1, we change the $1$-marked part $12$ to a $1$-marked part $\overline{12}$, and then we change the $2$-marked part $13$ to $2$-marked $12$. Then, we switch $1$-marked $\overline{15}$ to $1$-marked ${15}$ to get an overpartition in ${\widetilde{G}}_{7,6,5;2}(134)$.
\[\lambda=\begin{array}{cc}\left[
\begin{array}{cccccccccccccccc}
&2&&&5&&7&8&&&&12&&&\\
&2&&&5&6&&8&&&11&\mathbf{12}&&&\\
1&&&\bar{4}&&\bar{6}&&\bar{8}&&\overline{10}&&\mathbf{\overline{12}}&&&\mathbf{{15}}
\end{array}\right]&\begin{array}{c}3\\2\\1\end{array}
\end{array}.\]

Let $\lambda=(\lambda_1,\ldots,\lambda_m)\in{\widetilde{G}}_{N_1,\ldots,N_{k-1};i}$ be an overpartition of $n$. we consider the Gordon marking of $\lambda$, then there are $N_1$ $1$-marked parts in $G(\lambda)$, set $1\leq\lambda_{N_1}^{(1)}\leq\cdots\leq\lambda_1^{(1)}$. If $\lambda_p^{(1)}$ is an overlined part and $\lambda_{p-1}^{(1)}$ is a non-overlined part for $1<p<N_1$ or $\lambda_p^{(1)}$ is an overlined part for $p=1$, we can define the first dilation operation of $p$-th kind. We assume $\lambda_p^{(1)}=\overline{t}$, we consider the following two cases.

\noindent{\bf The first dilation operation of $p$-th kind:}

{\bf Case 1:} There are two parts of the same mark with underlying parts $t$ and $t-1$, we denote
this same mark by $r$. It should be noticed that there are no $1$-marked parts with underlying
part $t+1$ because of the choice of $\lambda_p^{(1)}$. We change $\lambda_p^{(1)}$ to a non-overlined part $t$ and replace
the $r$-marked part $t$ by an $r$-marked part $t+1$. Moreover, if there is a $1$-marked non-overlined part to
the right of $t$, then we switch it to an overlined part.

{\bf Case 2:} There are no two parts with underlying parts $t$ and $t-1$ that have the same mark.
We see that there is no $1$-marked part with underlying part $t+1$ because of the choice of
$\lambda_p^{(1)}$. We change $\lambda_p^{(1)}$ to a non-overlined part $t$ with mark $1$. We denote by $r$ the largest mark
of the parts equal to $t$, and replace the $r$-marked non-overlined part $t$ with an $r$-marked
non-overlined part $t+1$. Since $r$ is the largest mark of the parts equal to $t$, and $t+1$ is not
a $1$-marked part of $\lambda$, we see that $t+1$ cannot be a part with a mark not exceeding $r$. So we
may place the new part equal to $t$ in a position of mark $r$.
 Moreover, if there are $1$-marked parts larger than $t$, then we switch the non-overlined $1$-marked part next to $\lambda_p^{(1)}$ to an overlined part.

It has been proved in \cite{Chen-Sang-Shi-2013} that the first reduction operation of $p$-th kind and The first reduction dilation of $p$-th kind are mutually inverse mapping. What is more, it has been proved in \cite{Sang-Shi-2015} that such two operation preserve the three conditions of $\mathbb{\widetilde{G}}_{N_1,\ldots,N_{k-1};i}$.

Now we aim to give a proof of Theorem \ref{N-1}.

\noindent{\bf Proof of Theorem \ref{N-1}.}  Let $R_N$ denote the set of partitions with distinct negative parts which lay in $[-N,-1]$. To prove Theorem \ref{N-1}, we give a bijection between $\mathbb{\widetilde{G}}_{N_1,\ldots,N_{k-1};i}$ and $R_{N_1-1}\times\mathbb{\widetilde{B}}_{N_1,\ldots,N_{k-1};i}$.

On the one hand, for an overpartition $\lambda\in\mathbb{\widetilde{G}}_{N_1,\ldots,N_{k-1};i}(n)$, we construct a pair $(\tau,\mu)\in R_{N_1-1}\times\mathbb{\widetilde{B}}_{N_1,\ldots,N_{k-1};i}$ and $|\tau|+|\mu|=n$.

If there is no overlined part in $\lambda$, we set $\tau=\emptyset$ and $\mu=\lambda$.

If there exist overlined parts in $\lambda$, then the overlined parts must be marked $1$, set the overlined parts be $\lambda^{(1)}_{j_1}$, $\lambda^{(1)}_{j_2}$, \ldots, $\lambda^{(1)}_{j_s}$, where $1\leq j_1<j_2<\cdots<j_s<N_1$. Let $M_0=0$ and $M_l=\sum_{t=1}^lj_t$ where $1\leq l\leq s$.  Set $\lambda^{0}=\lambda$, we do the following operations for $l$ from $1$ to $s$ successively.

For each $l$, we iterate the following operation for $p$ from $1$ to $j_l$:  we do the first dilation of $p$-th kind of $\lambda^{M_{l-1}+p-1}$ and denote the resulting overpartition by $\lambda^{M_{l-1}+p}$.

After the $M_s$ operations above, we can get an overpartition $\lambda^{M_s}$ which contains no overlined part. Finally, set $\mu=\lambda^{M_s}$ and $\tau=(-j_1,-j_2,\ldots,-j_s)$, we get the desired pair $(\tau,\mu)$.

On the contrary, for a pair $(\tau,\mu)\in R_{N_1-1}\times\mathbb{\widetilde{B}}_{N_1,\ldots,N_{k-1};i}$ and $|\tau|+|\mu|=n$, we want to construct an overpartition $\lambda\in\mathbb{\widetilde{G}}_{N_1,\ldots,N_{k-1};i}(n)$ .

If $\tau=\emptyset$, we just set $\lambda=\mu$.

If $\tau\neq\emptyset$, then we set $\tau=(-j_1,-j_2,\ldots,-j_s)$ where $1\leq j_1<j_2<\cdots<j_s<N_1$. Let $M_l=\sum_{t=1}^lj_t$ where $1\leq l\leq s$ and set $\lambda^{M_s}=\lambda$, we do the following operations for $l$ from $s$ to $1$ successively.

For each $l$, we iterate the following operation for $p$ from $1$ to $j_l$:  we do the first reduction of $p$-th kind of $\mu^{M_{l}-p+1}$ and denote the resulting overpartition by $\mu^{M_{l}-p}$.

After the $M_s$ operations above, we can get an overpartition $\mu^{0}$ which belongs to $\mathbb{\widetilde{G}}_{N_1,\ldots,N_{k-1};i}(n)$. Finally, we just need to set $\lambda=\mu^{0}$.

Thus we complete the proof of  Theorem \ref{N-1}. \qed

Now we are in a position to prove Theorem \ref{GENERATING-nwe-2}. We first recall construction  called Gordon marking of a partition \cite{Kursungoz-2010}.

For a partition $\eta=(\eta_1,\eta_2,\ldots,\eta_m)$ where $1\leq\eta_1\leq \eta_2\leq\cdots\eta_m$. The Gordon marking of $\eta$ is an assignment of positive integers (marks) to $\eta$ from smallest part to largest part  such that the marks are as small as possible subject to equal or consecutive parts are assigned distinct marks. For example, the Gordon marking of $\eta=(1,1,2,2,2,3,4,5,5,6,6,6)$ is
\[(1_1,1_2,2_3,2_4,2_5,3_1,4_2,5_1,5_3,6_2,6_4,6_5),\]
which can   be represented by an array as follows, where column indicates the value of parts and the row (counted from bottom to top) indicates the mark.
\[
\begin{array}{ccc}\left[
\begin{array}{cccccccc}
&2&&&&6\\
&2&&&&6\\
&2&&&5&\\
1&&&4&&6\\
1&&3&&5&
\end{array}
\right]&\begin{array}{c}5\\4\\3\\2\\1\end{array}\end{array}
.\]

 Let $\mathbb{B}_{N_1,\ldots,N_{k-1};i}(n)$ denote the set of partitions $\eta$ counted by $B_{k,i}(m,n)$ that have $N_r$ $r$-marked parts in the Gordon marking of $\eta$ for $1\leq r\leq k-1$, where $m=\sum_{r=1}^{k-1}N_{r}$.

Define \[\mathbb{B}_{N_1,\ldots,N_{k-1};i}=\displaystyle\sum_{n\geq0}
\mathbb{B}_{N_1,\ldots,N_{k-1};i}(n).\]
 Applying Kur\c{s}ung\"{o}z's bijection for  Theorem \ref{R-R-G-left-o} in  \cite{Kursungoz-2010},  we see that for $k> i\geq1$,
\begin{equation}\label{GENREATING-EEK}
\sum_{\eta\in\mathbb{B}_{N_1,\ldots,N_{k-1};i}}
x^{\ell(\eta)}
q^{|\eta|}=
\frac{
q^{N^2_1+\cdots+N^2_{k-1}+N_{i}+\cdots+N_{k-1}}
x^{N_1+\cdots+N_{k-1}}}
{(q)_{N_1-N_2}\cdots(q)_{N_{k-2}-N_{k-1}}
(q^2;q^2)_{N_{k-1}}}.
\end{equation}

Notice that ordinary partitions can be regard as special overpartitions with no overlined part. For an  ordinary partition $\lambda$, it is easy to check that the Gordon marking for partition of $\lambda$ is the same as the Gordon marking for overpartition of $\lambda$. Thus we can get that $\widetilde{\mathbb{B}}_{N_1,\ldots,N_{k-1};i}=\mathbb{B}_{N_1,\ldots,N_{k-1};i}$. By \eqref{GENREATING-EEK}, we obtain that for $k> i\geq1$,
\begin{equation}\label{GENREATING-EEK-w}
\begin{split}
&\sum_{\mu\in\widetilde{\mathbb{B}}_{N_1,\ldots,N_{k-1};i}}
x^{\ell(\mu)}
q^{|\mu|}=\sum_{\eta\in\mathbb{B}_{N_1,\ldots,N_{k-1};i}}
x^{\ell(\eta)}
q^{|\eta|}\\
&=\frac{
q^{N^2_1+\cdots+N^2_{k-1}+N_{i}+\cdots+N_{k-1}}
x^{N_1+\cdots+N_{k-1}}}
{(q)_{N_1-N_2}\cdots(q)_{N_{k-2}-N_{k-1}}
(q^2;q^2)_{N_{k-1}}}
.\end{split}
\end{equation}

By Theorem \ref{N-1}, we get that for $k> i\geq1$,
\begin{equation}\label{N1w}
\begin{split}
&\sum_{\lambda\in\mathbb{\widetilde{G}}_{N_1,\ldots,N_{k-1};i}}x^{\ell(\lambda)}
q^{|\lambda|}=(-q^{1-N_1})_{N_1-1}\sum_{\mu\in
\mathbb{\widetilde{B}}_{N_1,\ldots,N_{k-1};i}}
x^{\ell(\mu)}
q^{|\mu|}\\
&=\frac{(-q^{1-N_1})_{N_1-1}
q^{N^2_1+\cdots+N^2_{k-1}+N_{i}+\cdots+N_{k-1}}
x^{N_1+\cdots+N_{k-1}}}
{(q)_{N_1-N_2}\cdots(q)_{N_{k-2}-N_{k-1}}
(q^2;q^2)_{N_{k-1}}}
.\end{split}
\end{equation}

Given the relation between ${\widetilde{G}}_{N_1,\ldots,N_{k-1};i}(n)$ and ${\widetilde{H}}_{N_1,\ldots,N_{k-1};i}(n)$ as stated in Lemma \ref{2Nnew1}, let $\mathbb{\widetilde{H}}_{N_1,\ldots,N_{k-1};i}=\cup_{n\geq0}\mathbb{\widetilde{H}}_{N_1,\ldots,N_{k-1};i}(n)$, we get the generating function for $\mathbb{\widetilde{H}}_{N_1,\ldots,N_{k-1};i}$.

\begin{thm}\label{GKI-1} For $k> i\geq1,$
\begin{equation}\label{GKI1}
\begin{split}
& \sum_{\pi\in\widetilde{\mathbb{H}}_{N_1,\ldots,N_{k-1};i}}
x^{\ell(\pi)}
q^{|\pi|}=\frac{(-q^{1-N_1})_{N_1-1}
q^{N^2_1+\cdots+N^2_{k-1}+N_{i+1}+\cdots+N_{k-1}}
x^{N_1+\cdots+N_{k-1}}}
{(q)_{N_1-N_2}\cdots(q)_{N_{k-2}-N_{k-1}}
(q^2;q^2)_{N_{k-1}}}.
\end{split}
\end{equation}
\end{thm}

\pf From the relation \eqref{2Nnew-3}, we deduce that for $1\leq i< k-1,$
\begin{equation}\label{GKI2}
\begin{split}
& \sum_{\pi\in\widetilde{\mathbb{H}}_{N_1,\ldots,N_{k-1};i}}
x^{\ell(\pi)}
q^{|\pi|}\\
&=\sum_{\lambda\in\widetilde{\mathbb{G}}_{N_1,\ldots,N_{k-1};i+1}}
x^{\ell(\lambda)}
q^{|\lambda|}\\
&=\frac{(-q^{1-N_1})_{N_1-1}
q^{N^2_1+\cdots+N^2_{k-1}+N_{i+1}+\cdots+N_{k-1}}
x^{N_1+\cdots+N_{k-1}}}
{(q)_{N_1-N_2}\cdots(q)_{N_{k-2}-N_{k-1}}
(q^2;q^2)_{N_{k-1}}}.
\end{split}
\end{equation}

For $i=k-1$, from \eqref{2Nnew-4} it follows that
\begin{equation}\label{GKI3}
\begin{split}
& \sum_{\pi\in\widetilde{\mathbb{H}}_{N_1,\ldots,N_{k-1};k-1}}
x^{\ell(\pi)}
q^{|\pi|}\\
&=\sum_{\lambda\in\widetilde{\mathbb{G}}_{N_1,\ldots,N_{k-1};1}}
x^{\ell(\lambda)}
q^{|\lambda|-N_1-\cdots-N_{k-1}}\\
&=\frac{(-q^{1-N_1})_{N_1-1}
q^{N^2_1+\cdots+N^2_{k-1}}
x^{N_1+\cdots+N_{k-1}}}
{(q)_{N_1-N_2}\cdots(q)_{N_{k-2}-N_{k-1}}
(q^2;q^2)_{N_{k-1}}}.
\end{split}
\end{equation}

Observe that the above formulas \eqref{GKI2} for $1\leq i< k-1$ and \eqref{GKI3} for $i=k-1$ take the same form \eqref{GKI1} as in the theorem. This completes the proof.  \qed

We are now ready to finish the proof of Theorem \ref{GENERATING-nwe-2}.

\noindent{\bf Proof of Theorem \ref{GENERATING-nwe-2}.} By the generating functions of $\widetilde{\mathbb{G}}_{N_1,\ldots,N_{k-1};i}$ and $\widetilde{\mathbb{H}}_{N_1,\ldots,N_{k-1};i}$ and relation \eqref{2N-1}, we find that

\begin{equation}\label{zuihou-1}
\begin{split}
&\quad \sum_{\delta\in\widetilde{\mathbb{E}}_{N_1,\ldots,N_{k-1};i}}
x^{\ell(\delta)}
q^{|\delta|}\\
&=\sum_{\lambda\in\widetilde{\mathbb{G}}_{N_1,\ldots,N_{k-1};i}}
x^{\ell(\lambda)}
q^{|\lambda|}+\sum_{\pi\in\widetilde{\mathbb{H}}_{N_1,\ldots,N_{k-1};i}}
x^{\ell(\pi)}
q^{|\pi|}\\
&=\frac{(-q^{1-N_1})_{N_1-1}
q^{N^2_1+\cdots+N^2_{k-1}+N_{i+1}+\cdots+N_{k-1}}(1+q^{N_i})x^{N_1+\cdots+N_{k-1}}}
{(q)_{N_1-N_2}\cdots(q)_{N_{k-2}-N_{k-1}}(q^{2};q^{2})_{N_{k-1}}}.
\end{split}
\end{equation}

Thus we have
\begin{equation}
\begin{split}
&\quad \displaystyle\sum_{m,n\geq0}\widetilde{E}_{k,i}(m,n)x^mq^n\\
&=\sum_{N_{1}\geq \cdots\geq N_{k-1}\geq 0}\sum_{\delta\in\widetilde{\mathbb{E}}_{N_1,\ldots,N_{k-1};i}}
x^{\ell(\delta)}
q^{|\delta|}\\
&=\sum_{N_{1}\geq \cdots\geq N_{k-1}\geq 0}\frac{(-q^{1-N_1})_{N_1-1}
q^{N^2_1+\cdots+N^2_{k-1}+N_{i+1}+\cdots+N_{k-1}}(1+q^{N_i})x^{N_1+\cdots+N_{k-1}}}
{(q)_{N_1-N_2}\cdots(q)_{N_{k-2}-N_{k-1}}(q^{2};q^{2})_{N_{k-1}}}.
\end{split}
\end{equation}

This completes the proof of Theorem \ref{GENERATING-nwe-2}.   \qed

\section{G\"ollnitz-Gordon marking}

In this section, we first recall the definition and related notations of G\"ollnitz-Gordon marking of an overpartition which are first introduced in \cite{He-Ji-Wang-Zhao-2016} and then we will give an outline of the proof of Theorem \ref{Gollnitz-even-e1}.

\begin{defi}[G\"ollnitz-Gordon marking]\label{pf-1}
For an overpartition $\lambda=(\lambda_1,\lambda_2,\cdots,\lambda_r)$ where $\lambda_1\leq\lambda_2\leq\cdots\leq\lambda_r$ also in the following order
\begin{equation}
\overline{1}<1<\overline{2}<2<\cdots.
\end{equation}
We assign the marks to parts also in the order above such that the marks are as small as possible subject to the following conditions:

$(1)$ All of non-overlined odd parts and overlined even parts are marked $1$;

$(2)$ The overlined odd part $\overline{2j+1}$ is assigned different mark with   $2j$  or $\overline{2j}$;

$(3)$ The mark of non-overlined even parts $2j+2$ is more complicated, we consider the following three cases: $($a$)$ The non-overlined even parts with equal  size  are assigned different marks; $($b$)$ If   $2j$ or $\overline{2j}$ or $\overline{2j+1}$ is marked by 1, or there do not exist $2j+1$ and $\overline{2j+2}$ in $\lambda$, then   $2j+2$ is assigned different mark with  $2j$, $\overline{2j}$ and $\overline{2j+1}$. $($c$)$ Otherwise,  $2j+2$ can not be assigned by $1$ and  be assigned different mark with the mark of $2j$ and $\overline{2j+1}$ except that the smallest mark assigned to $2j$ or $\overline{2j+1}$ can be used as the mark of $2j+2$.
\end{defi}

For example, we consider the overpartition \[\lambda=(1,1,\overline{2},2,\overline{3},
\overline{4},6,7,8,8,\overline{10},10,\overline{11},\overline{12},
\overline{13}).\]
The G\"ollnitz-Gordon marking of $\lambda$  is
\[GG(\lambda)=(1_1,1_1,{\overline{2}}_1,2_2,{\overline{3}}_3,\overline{4}_1,6_2,7_1,8_2,8_3,{\overline{10}}_1,{10}_2,{\overline{11}}_3,{\overline{12}}_1,{\overline{13}}_2).\]
Similar to the Gordon marking, we can represent G\"ollnitz-Gordon marking by an array where column indicates the value of parts, and the row (counted from bottom to top) indicates the mark, so the G\"ollnitz-Gordon marking of $\lambda$ can be expressed as follows
\begin{equation}\label{example1}
GG(\lambda)=\begin{array}{ccc}\left[
\begin{array}{cccccccccccccccccccccc}
 &&\overline{3}&&&&&8&&&\overline{11}&&\\
 &2&&&&6&&8&&10&&&\overline{13}\\
 1^2&\overline{2}&&\overline{4}&&&7&&&\overline{10}&&\overline{12}&
\end{array}
\right]&\begin{array}{c}3\\2\\1\end{array}\end{array}
\end{equation}
It is not hard to see that  the non-overlined odd parts and  overlined even parts in $\lambda$ have only  appeared in the first row of $GG(\lambda)$. Furthermore, the parts  in each row except for the first row are distinct. Moreover, there are only  odd parts repeated in the first row.   We use ${2j+1}^t$ to  denote that there are $t$ multiplies of $2j+1$  in the first row for $t\geq 2$, and ${\overline{2j+1}}^t$ to denote that there are one  $\overline{2j+1}$ and $t-1$ multiplies of   $2j+1$ in the first row.

Let $N_r$ be the number of  parts in the $r$-th row of $GG(\lambda)$. It is easy to find that $N_1\geq N_2\geq\cdots$. So we could define $n_r=N_r-N_{r+1}$ for any positive integer $r$.    For the example above, we have $N_1=7,\ N_2=5,\ N_3=3,\ n_1=N_1-N_2=2,\ n_2=N_2-N_3=2,\ n_3=N_3=3.$

For a part $\lambda_j$ of an overpartition $\lambda=(\lambda_1,
\lambda_2,\ldots, \lambda_r)$, we write $\lambda_j=\overline{a_j}$ to indicate that $\lambda_j$ is an overlined part and write $\lambda_j=a_j$ to indicate that $\lambda_j$ is a non-overlined part. Set $|\lambda_j|=a_j$.

We recall the cluster defined on the G\"ollnitz-Gordon marking  of an overpartition defined in \cite{He-Ji-Wang-Zhao-2016}.
\begin{defi}\label{pf-2}
For an overpartition $\lambda$, we denote the 1-marked parts in  the G\"ollnitz-Gordon marking $GG(\lambda)$ of  $\lambda$ by
\[\lambda^{(1)}_{N_1}\leq \lambda^{(1)}_{N_1-1}\leq \cdots \leq \lambda^{(1)}_{1}\]
also in the following order
\[\overline{1}<1<\overline{2}<2<\cdots.\]
  We proceed to decompose $\lambda$ into $N_1$ clusters according to 1-marked parts of $\lambda$ in the above order. In other word, we aim to define $N_1$-cluster corresponding to $\lambda^{(1)}_{N_1}$, $(N_1-1)$-cluster  corresponding to $\lambda^{(1)}_{N_1-1}$, $\cdots$, 1-cluster  corresponding to $\lambda^{(1)}_{1}$  consecutively. Denote  $j$-cluster by $\alpha^{(j)}$, then
  \[\lambda=\{\alpha^{(N_1)},\alpha^{(N_1-1)},\ldots, \alpha^{(1)}\}.\]
  We supposed that $(N_1+1)$-cluster is $\emptyset$.  The $j$-cluster $\alpha^{(j)}$ corresponding to $\lambda^{(1)}_{j}$ is defined by considering the following three cases.
   \begin{itemize}
\item[{\rm(1)}]If $|\lambda^{(1)}_{j}|=|\lambda^{(1)}_{j-1}|$, then $\alpha^{(j)}$ has only one part equal to $\lambda^{(1)}_{j}$.

\item[{\rm (2)}]If $\lambda^{(1)}_{j}$ is odd and $\mid \lambda^{(1)}_{j-1}\mid-\mid\lambda^{(1)}_{j}\mid=1$, then
$\alpha^{(j)}$ has only one part equal to $\lambda^{(1)}_{j}$.

\item[{\rm(3)}]Otherwise, $\alpha^{(j)}$ is a maximal length sub-overpartition $\alpha^{(j)}_1\leq\alpha^{(j)}_{2}\leq \cdots \leq \alpha^{(j)}_s $ where
     $\alpha^{(j)}_1=\lambda^{(1)}_{j}$ and for $2\leq b\leq s$, $\alpha^{(j)}_{b}$ is a $b$-marked part of $\lambda$ and not  in  $(j+1)$-cluster  satisfying the following conditions:

     \begin{itemize}
     \item[{\rm (i).}] If  $\alpha^{(j)}_{b-1}$ is odd, then $| \alpha^{(j)}_{b}|-| \alpha^{(j)}_{b-1}| =1$.
   \item[{\rm (ii).}] If  $|\alpha^{(j)}_{b-1}|$ is even and there does not exist $(b-1)$-marked part with size $|\alpha^{(j)}_{b-1}|+2$ in $\lambda$, then $0\leq|\alpha^{(j)}_{b}|-|\alpha^{(j)}_{b-1}|\leq 2$.
     \item[{\rm (iii).}]If  $\alpha^{(j)}_{b-1}$ is even and there exists a $(b-1)$-marked part with size $|\alpha^{(j)}_{b-1}|+2$ in $\lambda$, then $0\leq|\alpha^{(j)}_{b}|-|\alpha^{(j)}_{b-1}|\leq 1$.
     \end{itemize}
\end{itemize}

\end{defi}

For the example in \eqref{example1}, the overpartition $\lambda$ has seven clusters, namely
 \vspace{-5mm}
 \begin{figure}[H]
\centering
\label{new-17}
\includegraphics[width=12cm]{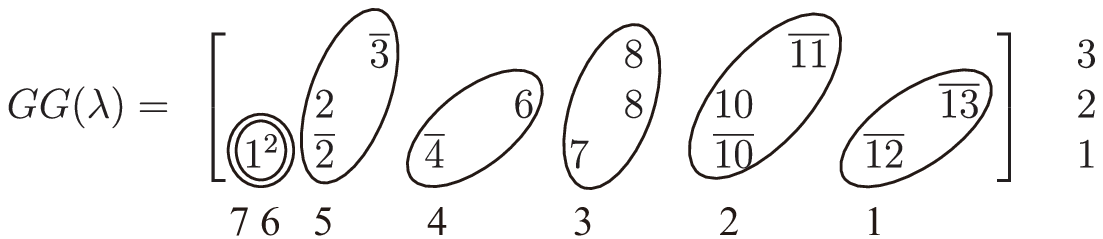}
\end{figure}
\vspace{-8mm}

From the definition of  cluster, it is easy to see that any overpartition $\lambda$ along with its G\"ollnitz-Gordon marking has a unique decomposition into non-overlapping clusters. The following proposition is useful, here we omit the proof for which you can see \cite{He-Ji-Wang-Zhao-2016}.

\begin{pro}\label{odd-part}
For an overpartition $\lambda$, if the G\"ollnitz-Gordon marking $GG(\lambda)$ is decomposed as follows  according to its clusters.
\[GG(\lambda)=\{\alpha^{(N_1)},\alpha^{(N_1-1)},\ldots, \alpha^{(1)}\},\] then for $1\leq t\leq N_1$, there exists at most one odd part in the $t$-cluster $\alpha^{(t)}$. Moreover, if there exists an odd part in $\alpha^{(t)}$, say $\alpha_{d}^{(t)}$, then $|\alpha_{d}^{(t)}|-|\alpha_1^{(t)}|\leq1$.
\end{pro}

The following proposition is also important.
\begin{pro}\label{xiangcha2}
For an overpartition $\lambda$, if the G\"ollnitz-Gordon marking $GG(\lambda)$ is decomposed as follows  according to its clusters.
\[GG(\lambda)=\{\alpha^{(N_1)},\alpha^{(N_1-1)},\ldots, \alpha^{(1)}\},\]
for $1\leq j< N_1$, we set
\[\alpha(j)=\{r\mid|\alpha_r^{(j)}|-|\alpha_r^{(j+1)}|\leq2 \text{ with restrict inequality if }\alpha_r^{(j)} \text{ or } \alpha_r^{(j+1)} \text{is odd}\},\]
then $\alpha(j)$ has at most one element.
\end{pro}

\pf Assume that $\alpha(j)$ has at least two elements.

For $r\in\alpha(j)$, we first prove that $\alpha_r^{(j)}$ and $\alpha_r^{(j+1)}$ are both even. Otherwise, $\alpha_r^{(j)}$ or $\alpha_r^{(j+1)}$ is odd, by the definition of G\"ollnitz-Gordon marking, we must have $r=1$. Now if $\alpha_1^{(j+1)}$ or $\alpha_1^{(j+1)}$ is odd, since $|\alpha_1^{(j)}|-|\alpha_1^{(j+1)}|<2$, by the definition of cluster, we get that $\alpha^{(j+1)}$  or $\alpha^{(j+1)}$ contains only one part, which contracts to the fact that $\alpha(j)$ has at least two elements. So $\alpha_r^{(j)}$ and $\alpha_r^{(j+1)}$ are both even.

Set $b$ and $c$ be the smallest element and the second smallest element in $\alpha(j)$, then $c>b$. Set $\alpha_c^{(j)}=2t+2$, then $\alpha_c^{(j+1)}=2t$. By the definition of G\"ollnitz-Gordon marking, we know that there is a $1$-marked $2t+1$ or $\overline{2t+2}$ and $c$ is the smallest mark of $\overline{2t}$, $2t$ and $\overline{2t+1}$. So there exist $2$-marked, $3$-marked,\ldots,$(c-1)$-marked $2t+2$ in $GG(\lambda)$. By the definition of cluster, it is easy to see that $\alpha_1^{(j)}=2t+1$ or $\overline{2t+2}$ and $\alpha_r^{(j)}=2t+2$ for $2\leq r\leq c$. Since $b<c$, so we have $\alpha_b^{(j)}=\overline{2t+2}$ or ${2t+2}$. But there is no $b$-marked $\overline{2t}$, $2t$ and $\overline{2t+1}$, which implies that $|\alpha_b^{(j)}|\leq2t-2$. So $|\alpha_b^{(j)}|-|\alpha_b^{(j+1)}|\geq4$, which contradicts to the choice of $b$.

Thus, the assumption is false, namely, $\alpha(j)$ has at most one element.  \qed

Let   $\mathbb{\widetilde{O}}_{N_1,\ldots,N_{k-1};i}(n)$ denote the set of overpartitions $\nu$ in $\mathbb{\widetilde{O}}_{k,i}(m,n)$ that have $N_r$ $r$-marked parts in their  G\"ollnitz-Gordon marking $GG(\nu)$ for $1\leq r\leq k-1$, where $\sum_{r=1}^{k-1}N_r=m$.

For an overpartition $\nu$, assume that the G\"ollnitz-Gordon marking representation of $\nu$ is decomposed as follows according to the clusters.
  \[GG(\nu)=\{\beta^{(N_1)},\beta^{(N_1-1)},\ldots, \beta^{(1)}\}.\]
 If there is an odd part in $\beta^{(p)}$ but no odd part in $\beta^{(p-1)}$ for $1< p\le N_1$ or there is an odd part in $\beta^{(p)}$ for $p=1$, we define the second dilation of $p$-th kind as follows.

{\noindent \bf The second dilation of $p$-th kind:}

We set the odd part in $\beta^{(p)}$ be $\beta_r^{(p)}$.

If $p=1$, we just change $\beta_r^{(1)}$ to an $r$-marked non-overlined (resp. overlined) even part with size $|\beta_r^{(1)}|+1$ if $\beta_r^{(1)}$ is overlined (resp. non-overlined).

If $1< p\le N_1$, we define
\[\beta(p-1)=\{j\mid|\alpha_j^{(p-1)}|-|\alpha_j^{(p)}|\leq2 \text{ with restrict inequality if }\alpha_j^{(p-1)} \text{ or } \alpha_j^{(p)} \text{is odd}\},\]
we consider following two cases:

{\bf Case 1:} $\beta(p-1)$ is empty, we first change $\beta_r^{(p)}$ to an $r$-marked non-overlined (resp. overlined) even part with size $|\beta_r^{(p)}|+1$ if $\beta_r^{(p)}$ is overlined (resp. non-overlined). Then we choose the part with the smallest size but the largest mark in $\beta^{(p-1)}$, set $b$-marked $\beta_b^{(p-1)}$, then we change $\beta_b^{(p-1)}$ to an $b$-marked non-overlined (resp. overlined) odd part with size $|\beta_b^{(p-1)}|+1$ if $\beta_b^{(p-1)}$ is overlined (resp. non-overlined).

{\bf Case 2:} $\beta(p-1)=\{b\}$, we first change $\beta_r^{(p)}$ to an $r$-marked non-overlined (resp. overlined) even part with size $|\beta_r^{(p)}|+1$ if $\beta_r^{(p)}$ is overlined (resp. non-overlined).Then we change $\beta_b^{(p-1)}$ to an $b$-marked non-overlined (resp. overlined) odd part with size $|\beta_b^{(p-1)}|+1$ if $\beta_b^{(p-1)}$ is overlined (resp. non-overlined).

As an example, let $GG(\nu)=\{\beta^{(7)},\beta^{(6)},\ldots, \beta^{(1)}\}$ be the following overpartition in $\mathbb{\widetilde{O}}_{7,5,3;3}(97)$.
 \vspace{-5mm}
 \begin{figure}[H]
\centering
\label{new-17}
\includegraphics[width=12cm]{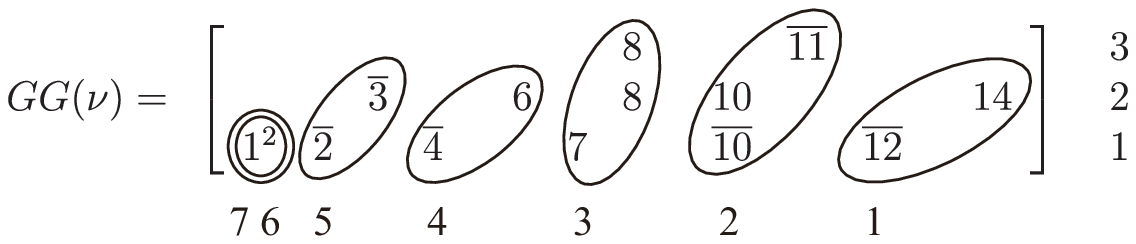}
\end{figure}
\vspace{-8mm}
There is an odd part $2$-marked $\overline{3}$ in $\beta^{(5)}$ but no odd part in $\beta^{(4)}$, so we can do the second dilation of $5$-th kind. Since $|\beta_1^{(4)}|-|\beta_1^{(5)}|=|\overline{4}|-|\overline{2}|=2$, we have $\beta(4)={1}$. So it belongs to Case 2. Now we first change $2$-marked $\overline{3}$ to a $2$-marked $4$, then we change $1$-marked $\overline{4}$ to a $1$-marked $5$. Thus we get the following overpartition belongs to $\mathbb{\widetilde{O}}_{7,5,3;3}(99)$.
\[
\begin{array}{ccc}\left[
\begin{array}{ccccccccccccccccccccccc}
 &&&&&&&8&&&\overline{11}&&\\
 &&&4&&6&&8&&10&&&&14\\
 1^2&\overline{2}&&&5&&7&&&\overline{10}&&\overline{12}&
\end{array}
\right]&\begin{array}{c}3\\2\\1\end{array}\end{array}
\]

For an overpartition $\nu$, assume that the G\"ollnitz-Gordon marking representation of $\nu$ is decomposed as follows according to the clusters.
  \[GG(\nu)=\{\beta^{(N_1)},\beta^{(N_1-1)},\ldots, \beta^{(1)}\}.\]
 If there is an odd part in $\beta^{(p-1)}$ but no odd part in $\beta^{(p)}$ for $1< p\le N_1$ or there is no odd part in $\beta^{(p)}$ for $p=1$, we define the second reduction of $p$-th kind as follows.

{\noindent \bf The second reduction of $p$-th kind:}

If $p=1$, we choose the part with the largest size but the smallest mark in $\beta^{(p)}$, set $b$-marked $\beta_b^{(p)}$, then we change $\beta_b^{(p)}$ to an $b$-marked non-overlined (resp. overlined) odd part with size $|\beta_b^{(p)}|-1$ if $\beta_b^{(p)}$ is overlined (resp. non-overlined).

If $1< p\le N_1$, set the odd part in $\beta^{(p-1)}$ be $\beta_r^{(p-1)}$ and define
\[\beta(p-1)=\{j\mid|\alpha_j^{(p-1)}|-|\alpha_j^{(p)}|\leq2 \text{ with restrict inequality if }\alpha_j^{(p-1)} \text{ or } \alpha_j^{(p)} \text{is odd}\},\]
we consider following two cases:

{\bf Case 1:} $\beta(p-1)$ is empty, we first change $\beta_r^{(p-1)}$ to an $r$-marked non-overlined (resp. overlined) even part with size $|\beta_r^{(p-1)}|-1$ if $\beta_r^{(p-1)}$ is overlined (resp. non-overlined). Then we choose the part with the smallest size but the largest mark in $\beta^{(p)}$, set $b$-marked $\beta_b^{(p)}$, then we change $\beta_b^{(p)}$ to an $b$-marked non-overlined (resp. overlined) odd part with size $|\beta_b^{(p)}|-1$ if $\beta_b^{(p)}$ is overlined (resp. non-overlined).

{\bf Case 2:} $\beta(p-1)=\{b\}$, we first change $\beta_r^{(p-1)}$ to an $r$-marked non-overlined (resp. overlined) even part with size $|\beta_r^{(p-1)}|-1$ if $\beta_r^{(p-1)}$ is overlined (resp. non-overlined). Then we change $\beta_b^{(p)}$ to an $b$-marked non-overlined (resp. overlined) odd part with size $|\beta_b^{(p)}|-1$ if $\beta_b^{(p)}$ is overlined (resp. non-overlined).

It is easy to check that the second dilation of $p$-th kind and  the second reduction of $p$-th kind are inverse mapping for each other. Now we want show that the second dilation of $p$-th kind and  the second reduction of $p$-th kind both preserve conditions $\rm{(1)}-\rm{(5)}$ in Theorem \ref{Gordon-Gollnitz-even-overp}. Such operations preserve conditions $\rm{(1)}-\rm{(3)}$ are easy to check. We first show that the second dilation preserves conditions $\rm{(4)}$ and $\rm{(5)}$ in Theorem \ref{Gordon-Gollnitz-even-overp}.

\begin{pro}
For an overpartition $\nu\in\mathbb{\widetilde{O}}_{N_1,\ldots,N_{k-1};i}(n)$, if the G\"ollnitz-Gordon marking representation of $\nu$ is decomposed as follows according to the clusters.
  \[GG(\nu)=\{\beta^{(N_1)},\beta^{(N_1-1)},\ldots, \beta^{(1)}\}.\]
 If there is an odd part $\beta_r^{(p)}$ in $\beta^{(p)}$ but no odd part in $\beta^{(p-1)}$ for $1< p\le N_1$ or there is an odd part in $\beta^{(p)}$ for $p=1$. Denote the overpartition after doing the second dilation of $p$-th kind of $\nu$ by $\nu'$, we shall show that if

\begin{equation}\label{condition4}
f_{2l}(\nu')+f_{\overline{2l}}(\nu')+f_{\overline{2l+1}}(\nu')+f_{2l+2}(\nu')= k-1
\end{equation}
for some $l$, then
\begin{equation}\label{condition4-1}
lf_{2l}(\nu')+lf_{\overline{2l}}(\nu')+lf_{\overline{2l+1}}(\nu')+(l+1)f_{2l+2}(\nu')\equiv V_{2l+1}(\nu')+i-1\pmod{2},
\end{equation}
and if
\begin{equation}\label{condition5}
 f_{2l+1}(\nu')\geq1\text{ and  }f_{2l+2}(\nu)= k-2
\end{equation}
for some $l$, then
\begin{equation}\label{condition5-1}
l+(l+1)f_{2l+2}(\nu')\equiv V_{2l+1}(\nu')+i-1\pmod{2}.
\end{equation}
\end{pro}

\pf Assume that $|\beta_r^{(p)}|=2t+1$.

If $p=1$, we change $r$-marked $\beta_r^{(1)}$ to an $r$-marked even part $2t+2$ (resp. $\overline{2t+2}$) if $\beta_r^{(1)}$ is $\overline{2t+1}$ (resp. $2t+1$).  For integer $l$ satisfying \eqref{condition4} or \eqref{condition5}, if $l\neq t$, it is easy to check that \eqref{condition4-1} or \eqref{condition5-1} is established.

If $\beta_r^{(1)}$ is $\overline{2t+1}$, it is easy to get $f_{2t}(\nu')=f_{2t}(\nu)$, $f_{\overline{2t}}(\nu')=f_{\overline{2t}}(\nu)$, $f_{\overline{2t+1}}(\nu')=f_{\overline{2t+1}}(\nu)-1$ , $f_{2t+2}(\nu')=f_{2t+2}(\nu)+1$ and $V_{2t+1}(\nu')=V_{2t+1}(\nu)-1$. So if \eqref{condition4} is established for $l=t$, we also have
\[f_{2t}(\nu)+f_{\overline{2t}}(\nu)+f_{\overline{2t+1}}(\nu)+f_{2t+2}(\nu)= k-1.\]
Thus,
\begin{equation*}
\begin{split}
&\quad tf_{2t}(\nu')+tf_{\overline{2t}}(\nu')+tf_{\overline{2t+1}}(\nu')+(t+1)f_{2t+2}(\nu')\\
&=tf_{2t}(\nu)+tf_{\overline{2t}}(\nu)+t(f_{\overline{2t+1}}(\nu)-1)+(t+1)(f_{2t+2}(\nu)+1)\\
&\equiv V_{2t+1}(\nu)+i-1-t+(t+1)\pmod{2}\\
&\equiv V_{2t+1}(\nu')+i-1\pmod{2}.
\end{split}
\end{equation*}
So \eqref{condition4-1} is established.

If $\beta_r^{(1)}$ is ${2t+1}$, it is east to get $f_{2t}(\nu')=f_{2t}(\nu)$, $f_{\overline{2t}}(\nu')=f_{\overline{2t}}(\nu)$, $f_{\overline{2t+1}}(\nu')=f_{\overline{2t+1}}(\nu)$ , $f_{2t+2}(\nu')=f_{2t+2}(\nu)$ and $V_{2t+1}(\nu')=V_{2t+1}(\nu)$. So if \eqref{condition4} is established for $l=t$, we also have
\[f_{2t}(\nu)+f_{\overline{2t}}(\nu)+f_{\overline{2t+1}}(\nu)+f_{2t+2}(\nu)= k-1.\]
Thus,
\begin{equation*}
\begin{split}
&\quad tf_{2t}(\nu')+tf_{\overline{2t}}(\nu')+tf_{\overline{2t+1}}(\nu')+(t+1)f_{2t+2}(\nu')\\
&=tf_{2t}(\nu)+tf_{\overline{2t}}(\nu)+t(f_{\overline{2t+1}}(\nu))+(t+1)(f_{2t+2}(\nu))\\
&\equiv V_{2t+1}(\nu)+i-1\pmod{2}\\
&\equiv V_{2t+1}(\nu')+i-1\pmod{2}.
\end{split}
\end{equation*}
So \eqref{condition4-1} is established.

Assume that \eqref{condition5} set up for $l=t$, then $f_{2t+1}(\nu')\geq1$. By the definition of cluster, we must have $\beta_r^{(1)}=2t+1$. So, $f_{2t+2}(\nu')=f_{2t+2}(\nu)=k-2$ and $V_{2t+1}(\nu')=V_{2t+1}(\nu)$. Thus,
\[t+(t+1)f_{2l+2}(\nu')=t+(t+1)f_{2l+2}(\nu)\equiv V_{2l+1}(\nu)+i-1\equiv V_{2l+1}(\nu')+i-1\pmod{2},\]
which is \eqref{condition5-1} for $l=t$.

Now we consider the case $1<p\le N_1$. Assume that the part we choose to change in $\beta^{(p-1)}$ is $\beta_b^{(p-1)}$. Set $|\beta_b^{(p-1)}|=2s+2$, we assert that $t\leq s$. By the proposition \ref{odd-part}, we know that $0\leq|\beta_r^{(p)}|-|\beta_1^{(p)}|\leq1$, namely $2t\le|\beta_1^{(p)}|\le2t+1$. Since $\beta_1^{(p-1)}$ is even, by the definition of G\"ollnitz-Gordon marking, we have  $|\beta_1^{(p-1)}|-|\beta_1^{(p)}|\geq1$ with restrict inequality if $\beta_1^{(p)}$ is even. Thus we have $|\beta_1^{(p-1)}|\geq2t+2$. By the definition of cluster, we know that $2s+2=|\beta_b^{(p-1)}|\geq|\beta_1^{(p-1)}|\geq2t+2$, which leads to our assertion.

For integer $l$ satisfying \eqref{condition4} or \eqref{condition5}, if $l\neq t$ and $s+1$, it is easy to check that \eqref{condition4-1} or \eqref{condition5-1} is established.

For integer $l$ satisfying \eqref{condition4} or \eqref{condition5}, if $l=t$ or $s+1$. If $s>t$, similar to the proof of $p=1$, it is easy to check that \eqref{condition4-1} or \eqref{condition5-1} is established. We just consider the case $s=t$ and $l=t$ or $t+1$.

For $s=t$, we have $2t+2=|\beta_b^{(p-1)}|\geq|\beta_1^{(p-1)}|\geq2t+2$, so $|\beta_1^{(p-1)}|=2t+2$. Since $2t\leq|\beta_1^{(p)}|\leq2t+1$, we know that $\beta_1^{(p-1)}=\overline{2t+2}$ and $\beta(p-1)=1$. So the second dilation of $p$-th kind only has following two cases:

{\bf Case 1:} $\beta_r^{(p)}$ is $2t+1$, we first change $r$-marked $2t+1$ to an $r$-marked $\overline{2t+2}$, and then we change the part $\beta_1^{(p-1)}=\overline{2t+2}$ to a $1$-marked $2t+3$.
We get that the occurrences of $2t,\ \overline{2t+1},\ 2t+2,\ \overline{2t+2},\ \overline{2t+3}$ and $2t+4$ in $\nu'$ are the same as those in $\nu$. Notice that $V_{2t+1}(\nu')=V_{2t+1}(\nu)$ and $V_{2t+3}(\nu')=V_{2t+3}(\nu)$, it is obvious that if \eqref{condition4} or \eqref{condition5} holds for $l=t$ or $t+1$, then \eqref{condition4-1} or \eqref{condition5-1} also holds for $l=t$ or $t+1$.

{\bf Case 2:} $\beta_r^{(p)}$ is $\overline{2t+1}$, we first change $r$-marked $\overline{2t+1}$ to an $r$-marked ${2t+2}$, and then we change the part $\beta_1^{(p-1)}=\overline{2t+2}$ to a $1$-marked $2t+3$. In such case, we must have $f_{2t+1}(\nu)=0$, otherwise we must have $\beta_1^{(p-1)}=2t+1$. But $\beta_1^{(p-1)}=\overline{2t+2}$, which leads to a contradiction. So we have $f_{2t+1}(\nu)=0$. Now it is easy to check that $f_{\overline{2t}}(\nu')=f_{\overline{2t}}(\nu)$, $f_{2t}(\nu')=f_{2t}(\nu)$, $f_{\overline{2t+1}}(\nu')=f_{\overline{2t+1}}(\nu)-1$, $f_{{2t+1}}(\nu')=f_{{2t+1}}(\nu)=0$, $f_{2t+2}(\nu')=f_{2t+2}(\nu)+1$, $f_{\overline{2t+2}}(\nu')=f_{\overline{2t+2}}(\nu)-1$, $f_{\overline{2t+3}}(\nu')=f_{\overline{2t+3}}(\nu)$, and $f_{2t+4}(\nu')=f_{2t+4}(\nu)$. What is more, we have $V_{2t+1}(\nu')=V_{2t+1}(\nu)-1$ and $V_{2t+3}(\nu')=V_{2t+3}(\nu)-2$ . Next we show that \eqref{condition4} and \eqref{condition5} lead to \eqref{condition4-1} and \eqref{condition5-1} respectively for $l=t$ or $t+1$.

Since $f_{{2t+1}}(\nu')=0$, so \eqref{condition5} can not hold for $l=t$. If  \eqref{condition4} holds for $l=t$, then we also have
\[f_{2t}(\nu)+f_{\overline{2t}}(\nu)+f_{\overline{2t+1}}(\nu)+f_{2t+2}(\nu)= k-1.\]
Thus,
\begin{equation*}
\begin{split}
&\quad tf_{2t}(\nu')+tf_{\overline{2t}}(\nu')+tf_{\overline{2t+1}}(\nu')+(t+1)f_{2t+2}(\nu')\\
&=tf_{2t}(\nu)+tf_{\overline{2t}}(\nu)+t(f_{\overline{2t+1}}(\nu)-1)+(t+1)(f_{2t+2}(\nu)+1)\\
&\equiv V_{2t+1}(\nu)+i-1+1\pmod{2}\\
&\equiv V_{2t+1}(\nu')+i-1\pmod{2}.
\end{split}
\end{equation*}
So \eqref{condition4-1} is established for $l=t$.

If \eqref{condition5} holds for $l=t+1$. Since $f_{2t+4}(\nu)=f_{2t+4}(\nu')=k-2$ and $f_{\overline{2t+2}}(\nu)=1$, we have $f_{2t+2}(\nu)=f_{\overline{2t+3}}(\nu)=0$ and
\[f_{2t+2}(\nu)+f_{\overline{2t+2}}(\nu)+f_{\overline{2t+3}}(\nu)+f_{2t+4}(\nu)= k-1.\]
So we can easily get
\begin{equation*}
\begin{split}
&\quad (t+1)+(t+2)f_{\overline{2t+4}}(\nu')\\
&=(t+1)f_{\overline{2t+2}}(\nu)+(t+2)f_{\overline{2t+4}}(\nu)\\
&=(t+1)f_{2t+2}(\nu)+(t+1)f_{\overline{2t+2}}(\nu)+(t+1)f_{\overline{2t+3}}(\nu)+(t+2)f_{2t+4}(\nu)\\
&\equiv V_{2t+3}(\nu)+i-1\pmod{2}\\
&\equiv V_{2t+3}(\nu')+i-1\pmod{2},
\end{split}
\end{equation*}
which leads to \eqref{condition5-1} for $l=t+1$.

Finally if \eqref{condition4} holds for $l=t+1$, we must have
\[f_{2t+2}(\nu)+f_{\overline{2t+2}}(\nu)+f_{\overline{2t+3}}(\nu)+f_{2t+4}(\nu)= k-1.\]
Thus,
\begin{equation*}
\begin{split}
&\quad (t+1)f_{2t+2}(\nu')+(t+1)f_{\overline{2t+2}}(\nu')+(t+1)f_{\overline{2t+3}}(\nu')+(t+2)f_{2t+4}(\nu')\\
&=(t+1)(f_{2t+2}(\nu)+1)+(t+1)(f_{\overline{2t+2}}(\nu)-1)+(t+1)f_{\overline{2t+3}}(\nu)+(t+2)f_{2t+4}(\nu)\\
&\equiv V_{2t+3}(\nu)+i-1\pmod{2}\\
&\equiv V_{2t+3}(\nu')+i-1\pmod{2}.
\end{split}
\end{equation*}
So \eqref{condition4-1} is established for $l=t+1$.

Now we complete the proof. \qed

If the G\"ollnitz-Gordon markings of an overpartition $\lambda$ are not exceed $k-1$, and $\lambda$ satisfies condition \eqref{condition4} $($or \eqref{condition5}$)$ but does not satisfies condition \eqref{condition4-1} $($or \eqref{condition5-1}$)$ for some $l$. Assume that we can do the second dilation of $\lambda$ and denote the resulting overpartition by $\nu$. Using the same method above, one can easily get that $\nu$ satisfies condition \eqref{condition4} $($or \eqref{condition5}$)$ but does not satisfies condition \eqref{condition4-1} $($or \eqref{condition5-1}$)$ for some $t$, we omit the details here.

Notice that the second dilation of $p$-th kind and the second reduction of $p$-th kind are inverse mapping for each other, so the second reduction of $p$-th kind also preserves conditions $\rm{(4)}$ and $\rm{(5)}$ in Theorem \ref{Gordon-Gollnitz-even-overp}.

Thus, we have the following proposition.
\begin{pro}
For an overpartition $\nu\in\mathbb{\widetilde{O}}_{N_1,\ldots,N_{k-1};i}(n)$ and it satisfies the condition of the second dilation of $p$-th kind \rm{(}resp. the second reduction of $p$-th kind\rm{)}, if we do the such operation for $\nu$ and denote the resulting overpartition by $\nu'$, then $\nu'\in\mathbb{\widetilde{O}}_{N_1,\ldots,N_{k-1};i}(n+2)$ \rm{(}resp. $\nu'\in\mathbb{\widetilde{O}}_{N_1,\ldots,N_{k-1};i}(n-2)$\rm{)} for $1<p\leq N_1$ or $\nu'\in\mathbb{\widetilde{O}}_{N_1,\ldots,N_{k-1};i}(n+1)$ \rm{(}resp. $\nu'\in\mathbb{\widetilde{O}}_{N_1,\ldots,N_{k-1};i}(n-1)$\rm{)} for $p=1$.
\end{pro}

At the end of this section, we will show the outline proof of Theorem \ref{Gollnitz-even-e1}. To  prove  Theorem \ref{Gollnitz-even-e1}, it suffice to show the following generating function.
\begin{equation}\label{equalthm}
\begin{split}
& \sum_{\nu\in\widetilde{\mathbb{O}}_{N_1,\ldots,N_{k-1};i}}
x^{\ell(\nu)}
q^{|\nu|}\\
&=\frac{(-q^{2-2N_1};q^2)_{N_1-1}(-q^{1-2N_1};q^2)_{N_1}
q^{2(N^2_1+\cdots+N^2_{k-1}+N_{i+1}+\cdots+N_{k-1})}
(1+q^{2N_i})x^{N_1+\cdots+N_{k-1}}}
{(q^2;q^2)_{N_1-N_2}\cdots(q^2;q^2)_{N_{k-2}-N_{k-1}}(q^{4};q^{4})_{N_{k-1}}}.
\end{split}
\end{equation}

Let $\mathbb{\widetilde{W}}_{N_1,\ldots,N_{k-1};i}(n)$ denote the set of overpartitions in $\mathbb{\widetilde{O}}_{N_1,\ldots,N_{k-1};i}(n)$ without odd parts. Set
\[\mathbb{\widetilde{W}}_{N_1,\ldots,N_{k-1};i}=\cup_{n\geq0}\mathbb{\widetilde{W}}_{N_1,\ldots,N_{k-1};i}(n),\]
we will show the following relation in the next section by constructing a bijection in terms of the second dilation and the second reduction.
\begin{thm}\label{finalN-1}
For $N_1\geq N_2\geq \cdots \geq N_{k-1}\geq 0$, we have
\begin{equation}\label{finalN1}
\sum_{\nu\in\mathbb{\widetilde{O}}_{N_1,\ldots,N_{k-1};i}}x^{\ell(\nu)}
q^{|\nu|}=(-q^{1-2N_1};q^2)_{N_1}\sum_{\omega\in
\mathbb{\widetilde{W}}_{N_1,\ldots,N_{k-1};i}}
x^{\ell(\omega)}
q^{|\omega|},
\end{equation}
\end{thm}
which leads a proof to \eqref{equalthm}.

\section{Proof of Theorem \ref{Gollnitz-even-e1}}

We first give a proof of Theorem \ref{finalN-1}.

 \noindent{\bf Proof of Theorem \ref{finalN-1}.} Let $Q_N$ denote the set of partitions with distinct negative odd parts which lay in $[1-2N,-1]$, We give a bijection between $\mathbb{\widetilde{O}}_{N_1,\ldots,N_{k-1};i}$ and $Q_{N_1}\times\mathbb{\widetilde{W}}_{N_1,\ldots,N_{k-1};i}$.

On the one hand, for an overpartition $\nu\in\mathbb{\widetilde{O}}_{N_1,\ldots,N_{k-1};i}(n)$, we construct a pair $(\xi,\omega)\in Q_{N_1}\times\mathbb{\widetilde{W}}_{N_1,\ldots,N_{k-1};i}$ and $|\xi|+|\omega|=n$.

If there is no odd part in $\nu$, we set $\xi=\emptyset$ and $\omega=\nu$.

If there exist odd parts in $\nu$, we consider the G\"ollnitz-Gordon marking $GG(\nu)=\{\beta^{(N_1)},\beta^{(N_1-1)},\ldots, \beta^{(1)}\}$. Then by proposition \ref{odd-part}, the odd parts must in different clusters of $GG(\nu)$, set such clusters be $\beta^{(j_1)}$, $\beta^{(j_2)}$, \ldots, $\beta^{(j_s)}$, where $1\leq j_1<j_2<\cdots<j_s\leq N_1$. Let $M_0=0$ and $M_l=\sum_{t=1}^lj_t$ where $1\leq l\leq s$.  Set $\nu^{0}=\nu$, we do the following operations for $l$ from $1$ to $s$ successively.

For each $l$, we iterate the following operation for $p$ from $1$ to $j_l$:  we do the second dilation of $p$-th kind of $\nu^{M_{l-1}+p-1}$ and denote the resulting overpartition by $\nu^{M_{l-1}+p}$.

After the $M_s$  operations above, we can get an overpartition $\nu^{M_s}$ which contains no odd part. Finally, set $\omega=\nu^{M_s}$ and $\xi=(1-2j_1,1-2j_2,\ldots,1-2j_s)$, we get a desired pair $(\xi,\omega)$.

On the contrary, for a pair $(\xi,\omega)\in Q_{N_1}\times\mathbb{\widetilde{W}}_{N_1,\ldots,N_{k-1};i}$ and $|\xi|+|\omega|=n$, we want to construct an overpartition $\nu\in\mathbb{\widetilde{O}}_{N_1,\ldots,N_{k-1};i}(n)$ .

If $\xi=\emptyset$, we just set $\nu=\omega$.

If $\xi\neq\emptyset$, then we set $\xi=(1-2j_1,1-2j_2,\ldots,1-2j_s)$ where $1\leq j_1<j_2<\cdots<j_s\leq N_1$. Let $M_l=\sum_{t=1}^lj_t$ where $1\leq l\leq s$ and set $\omega^{M_s}=\omega$, we do the following operations for $l$ from $s$ to $1$ successively.

For each $l$, we iterate the following operation for $p$ from $1$ to $j_l$:  we do the second reduction of $p$-th kind of $\omega^{M_{l}-p+1}$ and denote the resulting overpartition by $\omega^{M_{l}-p}$.

After the $M_s$ operations above, we can get an overpartition $\omega^{0}$ which belongs to $\mathbb{\widetilde{O}}_{N_1,\ldots,N_{k-1};i}(n)$. Finally, we just need to set $\nu=\omega^{0}$.

Thus we complete the proof of  Theorem \ref{finalN-1}. \qed

Now we are in a position to prove Theorem \ref{Gollnitz-even-e1}.

\noindent{\bf Proof of Theorem \ref{Gollnitz-even-e1}} Using \eqref{zuihou-1}, we shall prove that the generating function of   overpartitions in $\mathbb{\widetilde{W}}_{N_1,\ldots,N_{k-1};i}$ can be stated as follows.

\begin{lem}\label{finallww} For $k>i\geq1$, we have
\begin{equation}\label{GENREATING-EEww}
\sum_{\omega\in\mathbb{\widetilde{W}}_{N_1,\ldots,N_{k-1};i}}
x^{\ell(\omega)}
q^{|\omega|}=
\frac{(-q^{2-2N_1;q^2})_{N_1-1}
q^{2(N^2_1+\cdots+N^2_{k-1}+N_{i}+\cdots+N_{k-1})}
x^{N_1+\cdots+N_{k-1}}}
{(q^2;q^2)_{N_1-N_2}\cdots(q^2;q^2)_{N_{k-2}-N_{k-1}}
(q^{4};q^{4})_{N_{k-1}}}.
\end{equation}
\end{lem}

\pf
 We will establish a bijection between $\mathbb{\widetilde{E}}_{N_1,\ldots,N_{k-1};i}(n)$ and $\mathbb{\widetilde{W}}_{N_1,\ldots,N_{k-1};i}(2n)$.  For an overpartition $\delta=(\delta_1,\delta_2,\ldots,\delta_\ell)\in\mathbb{\widetilde{E}}_{N_1,\ldots,N_{k-1};i}(n)$, where $\delta_1\leq \delta_2\leq \cdots  \leq \delta_\ell$. Double each part in $\delta$ and denote the resulting overpartition by $\omega$, then $\omega=(2\delta_1,2\delta_2,\ldots,2\delta_\ell)$. It is easy to see that $|\omega|=2|\delta|=2n$ and $f_{\overline{1}}(\omega)+f_2(\omega)=f_2(\omega)=f_{1}(\delta)\leq i-1$. Since all parts in $\omega$ are even, so $f_{2t+1}(\omega)=f_{\overline{2t+1}}(\omega)=0$ and $V_{2t+1}(\omega)=V_{2t}(\omega)=V_{t}(\delta)$. If $f_{2t}(\omega)+f_{\overline{2t}}(\omega)+f_{\overline{2t+1}}(\omega)+f_{2t+2}(\omega)=k-1$, so $f_t(\delta)+f_{\overline{t}}(\delta)+f_{{t+1}}(\delta)=f_{2t}(\omega)+f_{\overline{2t}}(\omega)+f_{2t+2}(\omega)=k-1$. Notice that $\delta\in\mathbb{\widetilde{E}}_{N_1,\ldots,N_{k-1};i}(n)$, we have
 \begin{equation*}
\begin{split}
&\quad tf_{2t}(\omega)+tf_{\overline{2t}}(\omega)+tf_{\overline{2t+1}}(\omega)+(t+1)f_{2t+2}(\omega)\\
&=tf_{t}(\delta)+tf_{\overline{t}}(\delta)+(t+1)f_{t+1}(\delta)\\
&\equiv V_{t}(\delta)+i-1\pmod{2}\\
&\equiv V_{2t+1}(\omega)+i-1\pmod{2}.
\end{split}
\end{equation*}

 Next, we consider the Gordon marking of $\delta$ and the G\"ollnitz-Gordon marking of $\omega$. Recall the  Gordon marking of $\delta$ and  the G\"ollnitz-Gordon marking of $\omega$, it  can be checked that the mark $\omega_j=2\delta_j$ in the G\"ollnitz-Gordon marking of $\omega$ are the same as the mark of $\delta_j$ in the Gordon marking of $\delta$ where $1\leq j\leq \ell$. Hence we show that $\omega\in\mathbb{\widetilde{W}}_{N_1,\ldots,N_{k-1};i}(2n)$ and the process is inversive.  Thus  we have constructed a bijection between $\mathbb{\widetilde{E}}_{N_1,\ldots,N_{k-1};i}(n)$ and $\mathbb{\widetilde{W}}_{N_1,\ldots,N_{k-1};i}(2n)$.
Finally by \eqref{zuihou-1}, we have
\begin{eqnarray}
& &\sum_{\omega\in\mathbb{\widetilde{W}}_{N_1,\ldots,N_{k-1};i}}
x^{\ell(\omega)}
q^{|\omega|}\nonumber\\&=&
\sum_{\delta\in\mathbb{\widetilde{E}}_{N_1,\ldots,N_{k-1};i}}
x^{\ell(\delta)}q^{2|\delta|}\nonumber\\
&=&
\frac{(-q^{2-2N_1;q^2})_{N_1-1}
q^{2(N^2_1+\cdots+N^2_{k-1}+N_{i}+\cdots+N_{k-1})}
x^{N_1+\cdots+N_{k-1}}}
{(q^2;q^2)_{N_1-N_2}\cdots(q^2;q^2)_{N_{k-2}-N_{k-1}}
(q^{4};q^{4})_{N_{k-1}}}.\nonumber
\end{eqnarray}
Thus we complete the proof of Lemma \ref{finallww}. \qed

Applying Lemma \ref{finallww} in Theorem \ref{finalN-1}, we obtain
\begin{equation}
\begin{split}
&\quad\sum_{\nu\in\mathbb{\widetilde{O}}_{N_1,\ldots,N_{k-1};i}}x^{\ell(\nu)}
q^{|\nu|}=(-q^{1-2N_1};q^2)_{N_1}\sum_{\omega\in
\mathbb{\widetilde{W}}_{N_1,\ldots,N_{k-1};i}}
x^{\ell(\omega)}
q^{|\omega|}\\
&=\frac{(-q^{2-2N_1};q^2)_{N_1-1}(-q^{1-2N_1};q^2)_{N_1}
q^{2(N^2_1+\cdots+N^2_{k-1}+N_{i+1}+\cdots+N_{k-1})}
(1+q^{2N_i})x^{N_1+\cdots+N_{k-1}}}
{(q^2;q^2)_{N_1-N_2}\cdots(q^2;q^2)_{N_{k-2}-N_{k-1}}(q^{4};q^{4})_{N_{k-1}}},
\end{split}
\end{equation}
which is \eqref{equalthm}.

Hence we establish the following generating function of $\widetilde{O}_{k,i}(m,n)$.
\begin{equation}
\begin{split}
&\quad \displaystyle\sum_{m,n\geq0}\widetilde{O}_{k,i}(m,n)x^mq^n\\
&=\sum_{N_{1}\geq \cdots\geq N_{k-1}\geq 0}\sum_{\nu\in\widetilde{\mathbb{O}}_{N_1,\ldots,N_{k-1};i}}
x^{\ell(\nu)}
q^{|\nu|}\\
&=\sum_{N_{1}\geq \cdots\geq N_{k-1}\geq 0}\frac{(-q^{2-2N_1};q^2)_{N_1-1}(-q^{1-2N_1};q^2)_{N_1}
q^{2(N^2_1+\cdots+N^2_{k-1}+N_{i+1}+\cdots+N_{k-1})}
(1+q^{2N_i})x^{N_1+\cdots+N_{k-1}}}
{(q^2;q^2)_{N_1-N_2}\cdots(q^2;q^2)_{N_{k-2}-N_{k-1}}(q^{4};q^{4})_{N_{k-1}}}.
\end{split}
\end{equation}

This completes the proof of Theorem \ref{Gollnitz-even-e1}.   \qed


\vskip 0.5cm



\begin{thebibliography}{99}

\setlength{\itemsep}{-.8mm}

\bibitem{Agarwal-Andrews-Bressoud} A.K. Agarwal, G.E. Andrews and D.M. Bressoud, The Bailey lattice, J. Ind. Math. Soc. 51 (1987) 57--73.

\bibitem{Andrews-1967} G.E. Andrews, A generalization of the G\"{o}llnitz-Gordon partition theorem, Proc. Amer. Math. Soc. 18 (1967) 945--952.

\bibitem{Andrews-1974} G.E. Andrews, An analytic generalization of the Rogers-Ramanujan identities for odd moduli, Proc. Nat. Acad. Sci. USA 71 (1974) 4082--4085.

\bibitem{Andrews-1976} G.E. Andrews, The Theory of partitions, Addison-Wesley Publishing Co., 1976.

\bibitem{Andrews-1984} G.E. Andrews, Multiple series Rogers-Ramanujan type identities, Pacific J. Math. 114 (1984) 267--283.

\bibitem{Andrews-1986} G.E. Andrews, q-Series: Their development and application in analysis, number theory, combinatorics, physics, and computer algebra, American Mathematical Soc, 1986.

\bibitem{Bailey-1949} W.N. Bailey, Identities of the Rogers-Ramanujan type, Proc. London Math. Soc. 50(2) (1949) 1--10.

\bibitem{Bressoud-1979} D.M. Bressoud, A generalization of the Rogers-Ramanujan identities for all moduli, J. Combn. Theory, Ser. A 27 (1979) 64--68.

\bibitem{Bressoud-1980} D.M. Bressoud, Analytic and combinatorial generalizations of Rogers-Ramanujan identities, Mem. Amer. Math. Soc. 24(227) (1980) 54pp.

\bibitem{Bressoud-Ismail-Stanton-2000} D.M. Bressoud, M. Ismail, and D. Stanton, Change of base in Bailey pairs, Ramanujan J. 4 (2000) 435--453.

\bibitem{Chen-Sang-Shi-2013} W.Y.C. Chen, D.D.M. Sang and D.Y.H. Shi, The Rogers-Ramanujan-Gordon thoerem for overpartitions, Proc. London Math. Soc. 106(3) (2013) 1371--1393.

\bibitem{Chen-Sang-Shi-2015} W.Y.C. Chen, D.D.M. Sang and D.Y.H. Shi, An analogue of Bressoud's theorem of Rogers-Ramanujan type, Ramanujan J. 36 (2015) 69--80.


\bibitem{Corteel-2007} S. Corteel and O. Mallet, Overpartitions, lattice paths and Rogers-Ramanujan identities, J. Combin.
Theory Ser. A. 114(8) (2007) 1407--1437.

\bibitem{Corteel-2008} S. Corteel, J. Lovejoy, and O. Mallet, An extension to overpartitions of the Rogers-Ramanujan identities for even moduli, J. Number Theory 128 (2008) 1602--1621.



\bibitem{Gollnitz-1960} H. G\"ollnitz, Einfache Partionen, Diplomarbeit W. S. Gottingen, 65, 1960.

\bibitem{Gollnitz-1967} H. G\"ollnitz, Partitionen mit differenzenbedingungen, J. reine angew. Math. 225 (1967) 154--190.

\bibitem{Gordon-1961} B. Gordon, A combinatorial generalization of the Rogers-Ramanujan identities, Amer. J. Math 83 (1961) 393--399.

\bibitem{Gordon-1962} B.Gordon, Some Ramanujan-like continued fractions. Abstracts of Short Communications. Int. 1962. Congr. of Math., Stockholm, pp. 29--30, 1962.

\bibitem{Gordon-1965} B. Gordon, Some continued fractions of the Rogers-Ramanujan type, Duke Math. J. 31 (1965) 741--748.

\bibitem{He-Ji-Wang-Zhao-2016} Thomas Y. He, Kathy Q. Ji, Allison Y.F. Wang and
  Alice X.H. Zhao, The Andrews-G\"ollnitz-Gordon Theorem for Overpartitions, arXiv:1612.04960.


\bibitem{Kursungoz-2009} K. Kur\c{s}ung\"{o}z, Parity considerations in Andrews-Gordon identities, and the $k$-marked Durfee symbols, PhD thesis, Penn Sate University, 2009.

\bibitem{Kursungoz-2010} K. Kur\c{s}ung\"{o}z, Parity considerations in Andrews-Gordon identities, European J. Combin. 31 (2010) 976--1000.

\bibitem{Lovejoy-2003} J. Lovejoy, Gordon's theorem for overpartitions, J. Combin. Theory, Ser. A. 103 (2003) 393--401.

\bibitem{Lovejoy-2004} J. Lovejoy, Overpartition theorems of the Rogers-Ramanujan type, J. London Math. Soc. 69 (2004)
562--574.

\bibitem{Lovejoy-2007} J. Lovejoy, Partitions and overpartitions with attached parts, Arch. Math. (Basel) 88 (2007) 316--322.

\bibitem{Lovejoy-2010} J. Lovejoy, Partitions with rounded occurrences and attached parts, Ramanujan J. 23 (2010) 307--313.

\bibitem{Paule-1985} P. Paule, On identities of the Rogers-Ramanujan type, J. Math. Anal. Appl. 107 (1985) 255--284.

\bibitem{Paule-1987} P. Paule, A note on Bailey's lemma, J. Combin. Theory Ser. A 44 (1987) 164--167.

\bibitem{Rogers-1894} L.J. Rogers, Second memoir on the expansion of certain infinite products, Proc. London Math. Soc. 25 (1894) 318--343.

\bibitem{Sang-Shi-2015} D.D.M. Sang and D.Y.H. Shi, An Andrews-Gordon type identity for overpartitions, Ramanujan J. 37 (2015) 653--679.



\bibitem{Slater-1952} L.J. Slater, A new proof of Rogers's transformations of infinite series, Proc. London Math. Soc. 2(1) (1951) 460-475.

\end{thebibliography}
\end{document}